\documentclass[reqno,11pt]{amsart}
\usepackage{graphicx}
\usepackage{amscd}
\usepackage{amssymb}
\usepackage[mathscr]{eucal}
\numberwithin{equation}{section}
\allowdisplaybreaks[1]

\title[Refined Error Estimates for the Riccati Equation]{Refined Error Estimates for
the Riccati Equation with Applications to the Angular Teukolsky Equation}

\author[F.\ Finster]{Felix Finster}
\address{Fakult\"at f\"ur Mathematik \\ Universit\"at Regensburg \\ D-93040 Regensburg \\ Germany}
\email{finster@ur.de}

\author[J.\ Smoller]{Joel Smoller \\ \\ July 2013}
\thanks{J.S.\ is supported in part by the National Science Foundation,
Grant No.\ DMS-1105189.}
\address{Mathematics Department \\ The University of Michigan \\ Ann Arbor, MI 48109, USA}
\email{smoller@umich.edu}

\newtheorem{Def}{Def.}[section]
\newtheorem{Thm}[Def]{Theorem}
\newtheorem{Prp}[Def]{Proposition}
\newtheorem{Lemma}[Def]{Lemma}

\newcommand{\Thanks}{\vspace*{.5em} \noindent \thanks}
\newcommand{\Proof}{\begin{proof}}
\newcommand{\QED}{\end{proof} \noindent}

\renewcommand{\H}{\mathscr{H}}

\newcommand{\C}{\mathbb{C}}
\newcommand{\R}{\mathbb{R}}

\newcommand{\Z}{\mathbb{Z}}

\newcommand{\N}{\mathbb{N}}

\newcommand{\beq}{\begin{equation}}
\newcommand{\eeq}{\end{equation}}

\newcommand{\umax}{u_{\mbox{\tiny{max}}}}

\renewcommand{\O}{\mathscr{O}}
\newcommand{\dbar}{\mathchar'26\mkern-12mu \partial}

\DeclareMathOperator{\re}{Re}
\DeclareMathOperator{\im}{Im}

\begin{document}
\maketitle

\begin{abstract}
We derive refined rigorous error estimates for approximate solutions of
Sturm-Liouville and Riccati equations with real or complex potentials.
The approximate solutions include WKB approximations, Airy and parabolic
cylinder functions, and certain Bessel functions.
Our estimates are applied to solutions of the angular Teukolsky equation
with a complex aspherical parameter in a rotating black hole Kerr geometry.
\end{abstract}

\tableofcontents

\section{Introduction}
The Teukolsky equation arises in the study of electromagnetic, 
gravitational and neutrino-field perturbations in the Kerr geometry describing
a rotating black hole (see~\cite{chandra, teukolsky}).
In this equation, the {\em{spin}} of the wave enters as a parameter~$s \in \{0, \frac{1}{2}, 1,
\frac{3}{2}, 2, \ldots\}$ (the case~$s=0$ reduces to the scalar wave equation).
The Teukolsky equation can be separated into radial and angular parts, giving rise to a system of
coupled ODEs. Here we shall analyze the angular equation, also referred to as the spin-weighted
spheroidal wave equation. It can be written as the eigenvalue equation
\beq \label{eigen}
{\mathcal{A}} \,\Psi = \lambda \Psi\:,
\eeq
where the spin-weighted spheroidal wave operator~${\mathcal{A}}$ is an elliptic operator
with smooth coefficients on the unit sphere~$S^2$. More specifically, choosing polar
coordinates~$\vartheta \in (0, \pi)$ and~$\varphi \in [0, 2 \pi)$, we have (see for
example~\cite{whiting})
\[ {\mathcal{A}} \Theta = \lambda \Theta \quad \text{with} \quad
{\mathcal{A}} = - \frac{\partial}{\partial \cos \vartheta} \,\sin^2 \vartheta \,\frac{\partial}
{\partial \cos \vartheta}
+ \frac{1}{\sin^2 \vartheta} \left( \Omega\, \sin^2 \vartheta + i \frac{\partial}{\partial \varphi}
- s \cos \vartheta \right)^2 . \]
Here~$\Omega \in \C$ is the {\em{aspherical parameter}}. In the special case~$\Omega=0$,
we obtain the spin-weighted Laplacian, whose eigenvalues and eigenfunctions can be given
explicitly~\cite{goldberg}. In the case~$s=0$ and~$\Omega \neq 0$, one gets
the spheroidal wave operator as studied in~\cite{flammer}. Setting~$\Omega=0$ and~$s=0$,
one simply obtains the Laplacian on the sphere.
We are mainly interested in the cases~$s=1$ of an electromagnetic field and~$s=2$ of a gravitational field.

As the spin-weighted spheroidal wave operator is axisymmetric, we can separate out the
$\varphi$-dependence with a plane wave ansatz,
\[ \Psi(\vartheta, \varphi) = e^{-i k \varphi}\: \Theta(\vartheta) \qquad \text{with~$k \in \Z$}\:. \]
Then ${\mathcal{A}}$ becomes the ordinary differential operator
\beq \label{swswo}
{\mathcal{A}} =  - \frac{\partial}{\partial \cos \vartheta} \,\sin^2 \vartheta \,\frac{\partial}
{\partial \cos \vartheta}
+ \frac{1}{\sin^2 \vartheta} \left( \Omega\, \sin^2 \vartheta + k - s \cos \vartheta \right)^2  \:.
\eeq
To analyze the eigenvalue equation~\eqref{eigen}, we consider this operator on the Hilbert
space~$\H = L^2((-1,1), d\cos \vartheta)$ with domain of
definition~${\mathscr{D}}(\mathcal{A}) = C^\infty_0((-1,1))$.
In this formulation, the spheroidal wave equation also applies in the case of half-integer spin (to
describe neutrino or Rarita-Schwinger fields), if~$k$ is chosen to be a half-integer.
Thus in what follows, we fix the parameters~$s$ and~$k$ such that
\[ 2s \in \N_0 \qquad \text{and} \qquad k-s \in \Z\:. \]

In most applications, the aspherical parameter~$\Omega$ is real. However, having contour
methods for the Teukolsky equation in mind (similar as worked out in~\cite{weq} for the scalar wave equation),
we must consider the case that~$\Omega$ is complex. This leads to the major difficulty that
the potential in~\eqref{swswo} also becomes complex, so that the angular Teukolsky operator is no longer
a symmetric operator. At least, it suffices to consider the case when~$|\Omega|$ is large, whereas
the imaginary part of~$\Omega$ is uniformly bounded, i.e.\
\beq \label{Obounds}
|\Omega| > C \qquad \text{and} \qquad |\im \Omega| < c
\eeq
for suitable constants~$C$ and~$c$.
We are aiming at deriving a spectral representation for this non-symmetric
angular Teukolsky operator, which will involve complex eigenvalues and possibly Jordan chains.
In order to derive this spectral representation, we must have detailed knowledge of the
solutions of the Sturm-Liouville equation~\eqref{eigen}.
Our strategy for getting this detailed information
is to first construct approximate solutions by ``glueing together'' suitable
WKB, Airy, Bessel and parabolic cylinder functions, and then to derive rigorous error estimates.
The required properties of the special functions were worked out in~\cite{special}.
Our error estimates are based on the invariant region techniques in~\cite{invariant}.
These techniques need to be refined considerably in order to be applicable to the angular
Teukolsky equation.
Since these refined error estimates can be applied in a much more general context,
we organize this paper by first developing the general methods and then
applying them to the angular Teukolsky equation.

This paper is part of our program which aims at proving linear stability of the Kerr black hole
under perturbations of general spin (for another approach towards this goal we refer to~\cite{andersson+blue}). The estimates derived in the present paper play an important
role in this program. Indeed, they are a preparation for proving a spectral resolution
for the angular Teukolsky equation involving projectors onto finite-dimensional invariant
subspaces~\cite{tspectral} (related results on the angular operator for scalar waves can
be found in~\cite{angular, dyatlov}).
The next step will be to derive an integral representation
for the time evolution operator of the full Teukolsky equation in the same spirit as carried out in the
Schwarzschild geometry in~\cite{schdecay}.

The paper is organized as follows. We begin the analysis
by transforming the angular Teukolsky equation into Sturm-Liouville form with a complex
potential (Section~\ref{secsl}).
We then develop invariant region estimates for a general potential (Section~\ref{secinvreg}).
We proceed by deriving WKB estimates (Section~\ref{sec4}),
and then applying them to the angular Teukolsky equation (Section~\ref{sec5}).
In Section~\ref{sec6} we derive error estimates for parabolic cylinder approximations.
These include estimates for Airy approximations as a special case.
Section~\ref{sec7} is devoted to the properties of Bessel function solutions of
Sturm-Liouville equations with singular potentials.
Finally, in Section~\ref{secpoles} we use these properties to analyze solutions of
the angular Teukolsky equation near the poles at~$\vartheta = 0$ and~$\pi$.

\section{A Sturm-Liouville Operator with a Complex Potential} \label{secsl}
In order to bring the operator~\eqref{swswo} to the standard Sturm-Liouville form, we
first write the operator in the variable~$u=\vartheta \in (0,\pi)$,
\[ {\mathcal{A}} =  -\frac{1}{\sin u}\: \frac{d}{d u}\: \sin u\: \frac{d}{d u}
+ \frac{1}{\sin^2 u} \left( \Omega\, \sin^2 u + k - s \cos u \right)^2  \:. \]
Introducing the function $Y$ by
\beq \label{Ydef}
Y = \sqrt{\sin u}\: \Theta \:,
\eeq
we get the eigenvalue equation
\[ B\, \phi \;=\; \lambda \:\phi \;, \]
where
\begin{align*}
B &= -\frac{1}{\sqrt{\sin u}}\: \frac{d}{d u}\: \sin u\:
\frac{d}{d u}\:\frac{1}{\sqrt{\sin u}}
+\frac{1}{\sin^2 u}(\Omega \sin^2 u + k - s \cos u)^{2} \\
&= -\frac{d^2}{d u^2} + \frac{1}{2}\: \frac{\cos^2 u}{\sin^2 u}
-\sqrt{\sin u} \: \left(\frac{1}{\sqrt{\sin u}} \right)''
+\frac{1}{\sin^2 u}(\Omega \sin^2 u + k - s \cos u)^{2} \\
&= -\frac{d^2}{d u^2} - \frac{1}{4}\: \frac{\cos^2 u}{\sin^2 u}
- \frac{1}{2} +\frac{1}{\sin^2 u}(\Omega \sin^2 u + k - s \cos u)^{2} \:.
\end{align*}
Thus $\phi$ satisfies the Sturm-Liouville equation
\begin{equation} \label{5ode}
\left( -\frac{d^2}{du^2} + V \right) \phi \;=\; 0
\end{equation}
with the potential~$V$ given by
\begin{equation} \label{Vdef}
V \;=\; \Omega^2\: \sin^2 u + \left(k^2 + s^2 - \frac{1}{4} \right) \frac{1}{\sin^2 u}
- 2 s \Omega \cos u - 2 s k\, \frac{\cos u}{\sin^2 u} \:-\: \mu
\end{equation}
and~$\mu$ is the constant
\[ \mu \;=\; \lambda - 2 \Omega k + s^2 + \frac{1}{4} \:. \]
The transformation~(\ref{Ydef}) from~$\Theta$ to~$Y$ becomes
a unitary transformation if the integration measure in the
corresponding Hilbert spaces is transformed from~$\sin u\:du$ to~$du$.
Thus the eigenvalue problem~(\ref{eigen}) on~$\H$ is equivalent
to~(\ref{5ode}) on the Hilbert space~$L^2((0,\pi), du)$.

\section{General Invariant Region Estimates for the Riccati Flow} \label{secinvreg}
\subsection{An Invariant Disk Estimate}
Our method for getting estimates for solutions of the Sturm-Liouville
equation~\eqref{5ode} is to use invariant region estimates for the corresponding Riccati
equation. We here outline and improve the methods introduced in~\cite{invariant}.
Clearly, the solution space of the linear second order equation~\eqref{5ode} is two-dimensional.
For two solutions~$\phi_1$ and~$\phi_2$, the Wronskian~$w(\phi_1, \phi_2)$ defined by
\[ w(\phi_1, \phi_2) = \phi_1 \phi_2' - \phi_1' \phi_2 \]
is a constant. Integrating this equation, we can express one solution in terms of the other, e.g.
\[ \phi_2(u) = \phi_1(u) \left( \int^u \frac{w}{\phi_1^2} + \text{const} \right) \:. \]
Thus from one solution one gets the general solution by integration and taking
linear combinations. With this in mind, it suffices to get estimates for a {\em{particular}}
solution~$\phi$ of the Sturm-Liouville equation, which we can choose at our convenience.

Setting
\[ y = \frac{\phi'}{\phi}\:, \]
the function~$y$ satisfies the Riccati equation
\begin{equation} \label{riccati}
y'\;=\; V-y^2\:.
\end{equation}
Considering~$u$ as a time variable, the Riccati equation can be regarded as
describing a flow in the complex plane, the so-called {\em{Riccati flow}}.
In order to estimate~$y$, we want to find an approximate solution~$m(u)$ together
with a radius~$R(u)$ such that no solution~$y$ of the Riccati equation
may leave the circles with radius~$R$ centered at~$m$. More precisely, we
want that the implication
\[ \big| y(u_0) - m(u_0) \big| \leq R(u_0) \quad \Longrightarrow \quad
\big| y(u_1) - m(u_1) \big| \leq R(u_1) \]
holds for all~$u_1 > u_0$ and~$u_0, u_1 \in I$. We say that these circles
are {\em{invariant}} under the Riccati flow.
Decomposing~$m$ into real and imaginary parts,
\beq \label{mdef}
m(u) = \alpha(u) + i \beta(u)\:,
\eeq
our strategy is to prescribe the real part~$\alpha$, whereas the imaginary part~$\beta$
will be determined from our estimates. Then the functions~$U$ and~$\sigma$ defined by
\begin{align}
U &= \re V - \alpha^2 - \alpha'  \label{Udef} \\
\sigma(u) &= \exp \left( \int^u 2 \alpha \right) , \label{sigmadef}
\end{align}
which depend only on the known functions~$V$ and~$\alpha$,
can be considered as given functions. Moreover, we introduce the so-called
{\em{determinator}}~${\mathfrak{D}}$ by
\beq \label{Ddef}
{\mathfrak{D}} =  2 \alpha U + \frac{U'}{2} + \beta \im V \:.
\eeq
In our setting of a complex potential, the determinator involves~$\beta$ and
will thus be known only after computing the circles.
The following Theorem is a special case of~\cite[Theorem~3.3]{invariant} (obtained
by choosing~$W \equiv U$).

\begin{Thm}[{\bf{Invariant disk estimate}}] \label{thm1}
Assume that for a given function~$\alpha \in C^1(I)$ one of the following conditions holds:
\begin{itemize}
\item[{\bf{(A)}}] Defining real functions~$R$ and~$\beta$ on~$I$ by
\begin{align}
(R-\beta)(u) &= - \frac{1}{\sigma} \int^u \sigma
\im V \label{Rmb} \\[.5em] 
(R+\beta)(u) &= \frac{U(u)}{(R-\beta)(u)} \:, \label{Rpb}
\end{align}
assume that the function~$R-\beta$ has no zeros, $R \geq 0$, and
\beq \label{bc1}
(R-\beta) \:{\mathfrak{D}} \geq 0\:.
\eeq
\item[{\bf{(B)}}] Defining real functions~$R$ and~$\beta$ on~$I$ by
\begin{align}
(R+\beta)(u) &= \frac{1}{\sigma} \int^u \sigma
\im V \label{Rmb2} \\[.5em]
(R-\beta)(u) &= \frac{U(u)}{(R+\beta)(u)}\:, \label{Rpb2}
\end{align}
assume that the function~$R+\beta$ has no zeros, $R \geq 0$, and
\beq \label{bc2}
(R+\beta) \: {\mathfrak{D}} \;\geq\; 0 \:.
\eeq
\end{itemize}
Then the circle centered at~$m(u)=\alpha+i \beta$ with radius~$R(u)$ is invariant on~$I$
under the Riccati flow~\eqref{riccati}.
\end{Thm}
If this theorem applies and if the initial conditions~$y(u_0)$ lie
inside the invariant circles, we have obtained an approximate solution~$m$, \eqref{mdef}, together with
the rigorous error bound
\[ \big| y(u) - m(u) \big| \leq R(u) \qquad \text{for all~$u \geq u_0$}\:. \]

In order to apply the above theorems, we need to prescribe the function~$\alpha$.
When using Theorem~\ref{thm1}, the freedom in choosing~$\alpha$ must be used to
suitably adjust the sign of the determinator.
One method for constructing~$\alpha$ is to modify the potential~$V$ to a new
potential~$\tilde{V}$ for which the Sturm-Liouville equation has an explicit solution~$\tilde{\phi}$,
\beq \label{SLtilde}
\Big( -\frac{d^2}{du^2} + \tilde{V} \Big) \,\tilde{\phi} = 0 \:.
\eeq
We let~$\tilde{y}:= \tilde{\phi}'/\tilde{\phi}$ be the corresponding Riccati solution,
\beq \label{RICtilde}
\tilde{y}' = \tilde{V} - \tilde{y}^2\:,
\eeq
 and
define~$\alpha$ as the real part of~$\tilde{y}$. Denoting the imaginary part of~$\tilde{y}$
by~$\tilde{\beta}$, we thus have
\beq \label{tydef}
\tilde{y} = \alpha + i \tilde{\beta}\:.
\eeq
Writing the real and imaginary parts of the Riccati equation in~\eqref{tydef} separately, we obtain
\beq \label{abt}
\alpha' = \re \tilde{V}- \alpha^2 + \tilde{\beta}^2 \:,\qquad
\tilde{\beta}' = \im \tilde{V}- 2 \alpha \tilde{\beta}\:.
\eeq
In this situation, the determinator and the invariant disk estimates can be written in a particularly convenient form,
as we now explain. First, integrating the real part of~$\tilde{y}$, we find that the
function~$\sigma$, \eqref{sigmadef}, can be chosen as
\beq \label{sigmarel}
\sigma(u) = \exp \Big( \int^u 2 \alpha \Big) = \exp \Big( 2 \re \int^u \frac{\tilde{\phi}'}
{\tilde{\phi}} \Big) = |\tilde{\phi}|^2 \:.
\eeq
Moreover, applying the first equation in~\eqref{abt} to~\eqref{Udef}, we get
\beq \label{Urel}
U = \re (V-\tilde{V}) - \tilde{\beta}^2 \:.
\eeq
Differentiating~\eqref{Urel} and using
the second equation in~\eqref{abt}, we obtain
\[ U' = \re (V-\tilde{V})' + 4 \alpha \tilde{\beta}^2 - 2 \tilde{\beta} \im \tilde{V} \:. \]
Substituting this equation together with~\eqref{Urel} into~\eqref{Ddef} gives
(cf.~\cite[Lemma~3.4]{invariant})
\beq
{\mathfrak{D}} = 2 \alpha \re (V-\tilde{V}) \:+\:
\frac{1}{2} \re (V-\tilde{V})' - \tilde{\beta} \im \tilde{V}
+ \beta \im V\:. \label{Ddef2}
\eeq

\subsection{The $T$-Method}
The main difficulty in applying Theorem~\ref{thm1} is that one must satisfy the inequalities~\eqref{bc1}
or~\eqref{bc2} by giving the determinator a specific sign.
In the case~$|\beta|>R$, we know that Theorem~\ref{thm1} applies no matter what the sign of the
determinator is, because either~\eqref{bc1} or~\eqref{bc2} is satisfied.
This suggests that by suitably combining the cases~{\bf{(A)}} and~{\bf{(B)}},
one should obtain an estimate which does not involve the sign of~${\mathfrak{D}}$.
The next theorem achieves this goal.
It is motivated by the method developed in~\cite[Lemma~4.1]{angular} in the case of real potentials.
The method works only under the assumption that the function~$U$ given by~\eqref{Udef}
or~\eqref{Urel} is negative.

\begin{Thm} \label{thm2}
Assume that~$U<0$. We define~$\beta$ and~$R$ by
\beq \label{bReq}
\beta = \frac{\sqrt{|U|}}{2} \left(T + \frac{1}{T} \right) , \qquad
R = \frac{\sqrt{|U|}}{2} \left(T - \frac{1}{T} \right)
\eeq
where~$T \geq 1$ is a real-valued function which satisfies the differential inequality
\beq \label{Tpin}
\frac{T'}{T} \geq \left| \frac{\mathfrak{D}}{U} \right| -\frac{\im V}{\sqrt{|U|}} \:\frac{T^2-1}{2\, T} \:.
\eeq
Then the circle centered at~$m(u)=\alpha(u)+i \beta(u)$ with radius~$R(u)$ is invariant
under the Riccati flow~\eqref{riccati}.
\end{Thm}
\Proof Making the ansatz~\eqref{bReq}
with a free function~$T \geq 1$, the equations~\eqref{Rpb} and~\eqref{Rpb2} hold automatically.
Moreover, we see that~$0 \leq R < \beta$, so that if~${\mathfrak{D}} \leq 0$ we can apply case~{\bf{(A)}},
whereas if~${\mathfrak{D}} > 0$ we are in case~{\bf{(B)}}.
From~\eqref{Ddef} and~\eqref{sigmadef}, we find that
\beq \label{DUeq}
\frac{\mathfrak{D}}{U} = 2 \alpha + \frac{U'}{2 U} - \frac{\im V}{|U|}\: \beta
= \frac{(\sigma\, \sqrt{|U|})'}{\sigma\, \sqrt{|U|}} - \frac{\im V}{2 \sqrt{|U|}}\:
\left( T + \frac{1}{T} \right).
\eeq
In case~{\bf{(A)}}, differentiating~\eqref{Rmb} gives the equation
\[ \left( -\frac{\sigma\, \sqrt{|U|}}{T} \right)' = - \sigma \im \sqrt{V} \:. \]
Solving for~$T'/T$ gives
\[ \frac{T'}{T} = \frac{(\sigma\, \sqrt{|U|})'}{\sigma\, \sqrt{|U|}}
- \frac{\im V}{\sqrt{|U|}}\:T \:. \]
Substituting~\eqref{DUeq} and using~\eqref{bReq}, we obtain
\beq \label{case1}\frac{T'}{T} = \frac{\mathfrak{D}}{U}
- \frac{\im V}{|U|}\:R \:.
\eeq
In case~{\bf{(B)}}, we obtain similarly
\[ \left( \sigma\, \sqrt{|U|} \:T \right)' = \sigma \im \sqrt{V} \]
and thus
\[ \frac{T'}{T} = - \frac{(\sigma\, \sqrt{|U|})'}{\sigma\, \sqrt{|U|}} + \frac{\im V}{\sqrt{|U|}}\:
\frac{1}{T} \:. \]
Again using~\eqref{DUeq} and~\eqref{bReq}, we obtain
\beq \label{case2}
\frac{T'}{T} = - \frac{\mathfrak{D}}{U} - \frac{\im V}{|U|}\:R \:.
\eeq
Using that the quotient~$\mathfrak{D}/U$ is positive in case~{\bf{(A)}}
and negative in case~{\bf{(B)}}, we can combine~\eqref{case1} and~\eqref{case2}
to the differential equation
\[ \frac{T'}{T} =  \left| \frac{\mathfrak{D}}{U} \right| - \frac{\im V}{|U|}\:R \:, \]
which now holds independent of the sign of the determinator.
Using~\eqref{bReq}, this equation can be written as
\[ \frac{T'}{T} = \left| \frac{\mathfrak{D}}{U} \right| -\frac{\im V}{\sqrt{|U|}} \:\frac{T^2-1}{2\,T} \:. \]
If~$T$ solves this equation, then we know from Theorem~\ref{thm1} that we have invariant circles for the
Riccati flow. Replacing the equality by an inequality, the function~$T$ grows faster.
Since increasing~$T$ increases the circle defined by~\eqref{bReq}, we again obtain
invariant regions.
\QED

The next theorem gives a convenient method for constructing a solution
of the inequality~\eqref{Tpin}.
\begin{Thm} \label{thmT}
Assume that~$U<0$. We choose a real-valued function~$g$
and define the function~$T$ by
\[ \log T(u) = \int^u  E\:, \]
where
\[ E = \big| E_1 + E_2 + E_3 \big| + E_4 \]
and
\begin{align*}
E_1 &:= \frac{1}{2\, |U|} \left(  4 \alpha \re (V-\tilde{V}) + \re (V-\tilde{V})' \right) \\
E_2 &:= \frac{\tilde{\beta}}{|U|}\, \im (V-\tilde{V}) \\
E_3 &:= -\frac{\im V}{|U|}\: \frac{\re(V-\tilde{V})}{\sqrt{|U|} + \tilde{\beta}} \\
E_4 &:= \frac{|\im V|}{\sqrt{|U|}} \: g(u) \:. 
\end{align*}
Then the circle centered at~$m(u)=\alpha(u)+i \beta(u)$ with radius~$R(u)$ is invariant
under the Riccati flow~\eqref{riccati}, provided that the following condition holds:
\beq \label{gdef}
\left\{ \begin{array}{ll} {\displaystyle g \geq -\frac{T-1}{T}} &  \text{if~$\im V \geq 0$} \\[.8em]
g \geq T-1 \quad & \text{if~$\im V < 0$\:.} \end{array} \right.
\eeq
\end{Thm}
\Proof According to the first equation in~\eqref{bReq},
\[ \left| \beta - \sqrt{|U|} \right| = \sqrt{|U|} \:\frac{(T-1)^2}{2\,T} \:. \]
Using this equation in~\eqref{Ddef2}, we obtain
\begin{align*}
|{\mathfrak{D}}| &\leq  \Big| 2 \alpha \re (V-\tilde{V}) \:+\:
\frac{1}{2}\, \re (V-\tilde{V})' \\
&\qquad + \tilde{\beta} \im (V-\tilde{V})
+ \left( \sqrt{|U|} - \tilde{\beta} \right) \im V \Big|+ \sqrt{|U|} \,|\im V|\:\frac{(T-1)^2}{2\,T}\:.
\end{align*}
Applying the identities
\[ \sqrt{|U|} - \tilde{\beta} = \frac{|U|^2 - \tilde{\beta}^2}{\sqrt{|U|} + \tilde{\beta}}
= -\frac{\re(V-\tilde{V})}{\sqrt{|U|} + \tilde{\beta}} \]
(where in the last step we applied~\eqref{Urel} and used that~$U<0$),
the right side of~\eqref{Tpin} can be estimated by
\begin{align*}
\left| \frac{\mathfrak{D}}{U} \right|-\frac{\im V}{\sqrt{|U|}} \:\frac{T^2-1}{2\, T} 
\leq |E_1+E_2+E_3| + \frac{|\im V|}{\sqrt{|U|}}\:\frac{(T-1)^2}{2\,T}
-\frac{\im V}{\sqrt{|U|}} \:\frac{T^2-1}{2\, T} \:.
\end{align*}
Simplifying the last two summands in the two cases~$\im V \geq 0$ and~$\im V < 0$ gives the result.
\QED

\subsection{The $\kappa$-Method} \label{seckappa}
We now explain an alternative method for getting invariant region estimates.
This method is designed for the case when~$|\beta|<R$.
In this case, the factors~$R \mp \beta$ in~\eqref{bc1} and~\eqref{bc2} have the same sign.
Therefore, Theorem~\ref{thm1} applies only if the determinator has has the right sign.
In order to arrange the correct sign of the determinator, we must work with
driving functions (for details see Section~\ref{sec42}). When doing this, we know a-priori whether
we want to apply Theorem~\ref{thm1} in case~{\bf{(A)}} or~{\bf{(B)}}.
With this in mind, we may now restrict attention to a fixed case~{\bf{(A)}} or~{\bf{(B)}}.
In order to treat both cases at once, whenever
we use the symbols~$\pm$ or~$\mp$, the upper and lower signs refer to the
cases~{\bf{(A)}} and~${\bf{(B)}}$, respectively.
Differentiating~\eqref{Rmb} and~\eqref{Rmb2} and
using the form of~$\sigma$, \eqref{sigmadef}, we obtain
\[ (R \mp \beta)' = -2 \alpha\, (R \mp \beta) \mp \im V \:. \]
Combining this differential equation with the second equation in~\eqref{abt}, we get
\[ \big( \beta \mp R - \tilde{\beta} \big)' = \im (V-\tilde{V})
-2 \alpha\, \big( \beta \mp R - \tilde{\beta} \big) . \]
This differential equation can be integrated. Again using~\eqref{sigmadef}, we
obtain
\begin{gather}
\beta \mp R - \tilde{\beta} = \kappa \qquad \text{with} \label{tbeq} \\
\kappa := \frac{1}{\sigma} \left( \int^u \sigma \im (V-\tilde{V}) + C \right) , \label{tbeq2}
\end{gather}
where the integration constant~$C$ must be chosen such that~\eqref{tbeq} holds initially.
Solving~\eqref{tbeq} for~$\beta$ and using the resulting equation in~\eqref{Ddef2}
gives
\beq
{\mathfrak{D}} = 2 \alpha\, \re (V-\tilde{V})
+ \frac{1}{2}\, \re (V-\tilde{V})'
+ \tilde{\beta} \im (V-\tilde{V})
+ (\kappa \pm R) \im V \:. \label{Ddef3}
\eeq
The combination~$\kappa \pm R$ in~\eqref{Ddef3} has the following useful representation.
\begin{Lemma} \label{lemmaRpmb}
The function~$\kappa \pm R$ is given by
\beq \label{Rpmb}
\kappa \pm R = \frac{\kappa^2 - \re (V-\tilde{V})}{2 \,(\tilde{\beta} + \kappa)} \:.
\eeq
\end{Lemma}
\Proof
According to~\eqref{tbeq} and~\eqref{Rpb}, \eqref{Rpb2},
\[ R \mp \beta = \mp (\tilde{\beta} + \kappa) \:,\qquad
R \pm \beta = \frac{U}{R \mp \beta} = \mp \frac{U}{\tilde{\beta} + \kappa} \]
and thus
\[ R =  \mp \frac{1}{2} \left( (\tilde{\beta} + \kappa) + \frac{U}{\tilde{\beta} + \kappa} \right) 
= \mp \frac{U + (\tilde{\beta} + \kappa)^2}{2 \,(\tilde{\beta} + \kappa)} \:. \]
It follows that
\[ \kappa \pm R = \kappa - \frac{U + (\tilde{\beta} + \kappa)^2}{2 \,(\tilde{\beta} + \kappa)}
= \frac{2 \kappa \tilde{\beta} + 2 \kappa^2 -U - (\tilde{\beta} + \kappa)^2}{2 \,(\tilde{\beta} + \kappa)}
= \frac{\kappa^2 -U - \tilde{\beta}^2}{2 \,(\tilde{\beta} + \kappa)} , \]
and using~\eqref{Urel} gives the result.
\QED
The above relations give the following method for getting invariant region estimates.
First, we choose an approximate potential~$\tilde{V}$ having an explicit solution~$\tilde{y}=
\alpha + i \tilde{\beta}$. Next, we compute~$\sigma$ by~\eqref{sigmadef} or~\eqref{sigmarel}
and compute the integral~\eqref{tbeq2} to obtain $\kappa$.
The identity~\eqref{Rpmb} gives the quantity~$\kappa \pm R$.
Substituting this result into~\eqref{Ddef3}, we get an explicit formula for the
determinator. Instead of explicit computations, one can clearly 
work with inequalities to obtain estimates of the determinator.
The key point is to use the freedom in choosing~$\tilde{V}$ to give the
determinator a definite sign. Once this has been accomplished, we
can apply Theorem~\ref{thm1} in cases~{\bf{(A)}} or~{\bf{(B)}}.

The method so far has the disadvantage that the function~$\tilde{\beta}+\kappa$
in the denominator in~\eqref{Rpmb} may become small, in which case the
summand~$(\kappa \pm R) \im V$ in the determinator~\eqref{Ddef3} gets out of control.
Our method for avoiding this problem is to increase~$\kappa$ in such a way that the solution
stays inside the resulting disk. This method only works in case~${\bf{(B)}}$
of Theorem~\ref{thm1}.
\begin{Prp} \label{prpnewkappa}
Assume that~$y$ is a solution of the Riccati equation~\eqref{riccati}
in the upper half plane~$\im y > 0$. Moreover, assume that~${\mathfrak{D}}>0$.
For an increasing function~$g$ we set
\beq \label{kset}
\kappa(u) = \frac{g(u)}{\sigma(u)} + \frac{1}{\sigma} \int^u \sigma \im(V-\tilde{V})
\eeq
and choose~$R$ and~$\beta$ according to~\eqref{tbeq} and~\eqref{Rpb2},
\beq \label{ninvdisk}
R+\beta = \tilde{\beta}+\kappa \:,\qquad R-\beta = \frac{U}{R+\beta}\:.
\eeq
Then the circle centered at~$m=\alpha+i \beta$ with radius~$R$ is invariant on~$I$
under the Riccati flow. Moreover, Lemma~\ref{lemmaRpmb} remains valid.
\end{Prp}
\Proof According to Theorem~\ref{thm1}~{\bf{(B)}} and~\eqref{tbeq},
the identities~\eqref{ninvdisk} give rise to invariant disk estimates if we choose
\beq \label{relB}
\tilde{\beta} + \kappa = \tilde{\beta} + \frac{\text{const}}{\sigma(u)} + \frac{1}{\sigma}
\int^u \sigma \im(V-\tilde{V}) \:.
\eeq
If the constant is increased, the upper point~$R+\beta$ of the circle moves up.
In the case~$\beta-R \geq 0$, the second equation in~\eqref{ninvdisk} implies
that the lower point~$\beta-R$ of the circle moves down. As a consequence,
the disk increases if the constant is made larger.
Likewise, in the case~$\beta-R < 0$, the circle intersects the axis~$\im y=0$ in the two
points~$\alpha \pm \sqrt{U}$, which do not change if the constant is increased.
As a consequence, the intersection of the disk with the upper half plane increases
if the constant is made larger. Thus in both cases, the solution~$y(u)$ stays inside the
disk if the constant is increased.

We next subdivide the interval~$I$ into
subintervals. On each subinterval, we may use the formula~\eqref{relB} with an increasing sequence
of constants. Letting the number of subintervals tend to infinity, we conclude
that we obtain an invariant region estimate if the constant in~\eqref{relB} is replaced by a
monotone increasing function~$g(u)$.
\QED

\subsection{Lower Bounds for~$\im y$}
We begin with an estimate in the case when~$\im V$ is positive.
\begin{Lemma} \label{lemma36}
Suppose that~$y$ is a solution of the Riccati equation~\eqref{riccati} for a potential with the property
\beq \label{imVp}
\im V > 0\:.
\eeq
Assume furthermore that~$\im y(u_0) > 0$. Then
\beq
\im y(u) \geq \im y(u_0)\: \exp \left( -2 \int_{u_0}^u \re y \right) . \label{yes1}
\eeq
Moreover, the Riccati flow preserves the inequality
\beq
\im y(u) \geq \inf_{[u_0, u]} \frac{\im V}{2\, \re y}\:. \label{yes2}
\eeq
\end{Lemma}
\Proof Taking the imaginary part of~\eqref{riccati} gives
\beq
\im y' = \im V - 2 \re y\, \im y \:. \label{imric}
\eeq
From~\eqref{imVp}, we obtain
\[ \log' |\im y| \geq -2 \re y \:. \]
Integration gives~\eqref{yes1}. In particular, $\im y$ stays positive.

For the proof of~\eqref{yes2} we assume conversely that this inequality holds at
some~$u_1>u_0$ but is violated for some~$u_2 > u_1$.
Thus, denoting the difference of the left and right side of~\eqref{yes2} by~$g$,
we know that~$g(u_1) \geq 0$ and~$g(u_2)<0$. By continuity, there is
a largest number~$\bar{u} \in [u_1, u_2)$ with~$g(\bar{u})=0$.
According to the mean value theorem, there is~$v \in [\bar{u}, u_2]$ with~$g'(v) =g(u_2)/(u_2-\bar{u})<0$.
Since the function on the right is monotone decreasing in~$u$,
this implies that~$\im y'(v) < 0$.
Using~\eqref{imric}, we obtain at~$v$
\[ 0 > \im y'(v) \ =  \im V(v) - 2 \re y(v)\, \im y(v)\:. \]
If~$\re y \leq 0$, the infimum in~\eqref{yes2} is also negative, so that there is nothing to prove.
In the remaining case~$\re y>0$, we can solve for~$\im y$ to obtain
\[ \im y(v) > \frac{\im V(v)}{2 \re y(v)}  \geq \inf_{[u_0, v]} \frac{\im V}{2 \re y} \:. \]
Hence~$g(v)>0$, a contradiction.
\QED

The following estimate applies even in the case when~$\im V$ is negative. The method is to combine
a Gr\"onwall estimate with a differential equation for~$\im y$.
\begin{Lemma} \label{lemmagronwall} Let~$y$ be a solution of the Riccati equation~\eqref{riccati}
on an interval~$[u_-, u_+]$ and
\[ \max_{[u_-, u_+]} \sqrt{|V|}\: (u_+-u_-) \leq c\:. \]
Assume that~$\im y(u_-) \geq 0$. Then there is a constant~$C$ depending only on~$c$ such that
\[ \im y(u) \geq \frac{1}{C}\, \im {y}(u_-) - C\: (u-u_-) \:\Big| \min_{[u_-, u]} \im V \Big| \:. \]
\end{Lemma}
\Proof Let~$\phi(u) = \exp(\int^u y)$ be the corresponding solution of the
Sturm-Liouville equation~\eqref{5ode}.
Setting~$\kappa =\max_{[u_-, u_+]} |V|^\frac{1}{2}$, we write the Sturm-Liouville equation as the first
order system
\[ \Psi'(u) = \begin{pmatrix} 0 & \kappa \\ V/\kappa & 0 \end{pmatrix} \Psi(u)
\qquad \text{with} \qquad \Psi(u) := \begin{pmatrix} \kappa \,{\phi}(u) \\ {\phi}'(u) \end{pmatrix} . \]
Using that
\[ \int_{u_-}^{u_+} \left( \kappa + \frac{|V|}{\kappa} \right) du \leq \max_{[u_-, u_+]} \sqrt{|V|}
\:(u_+-u_-) \leq c \:, \]
a Gr\"onwall estimate yields
\beq \label{gronwall}
\frac{1}{c_2} \: \|\Psi(u_-)\| \leq \|\Psi(u) \| \leq c_2 \,\|\Psi(u_-)\| \:,
\eeq
where~$c_2$ depends only on~$c$.
This inequality bounds the combination~$\kappa^2 |{\phi}|^2 + |{\phi}'|^2$
from above and below. However, it does not rule out zeros of the function~${\phi}$.
To this end, we differentiate the identity
\[ \im (\overline{{\phi}} \,{\phi}' ) = \im (|{\phi}|^2 \, {y} ) 
= |{\phi}|^2 \im {y} \]
to obtain the differential equation
\[ \frac{d}{du} \left( |{\phi}|^2 \im {y} \right) = \im V\, |{\phi}|^2 \:. \]
Integrating this differential equation, we obtain
\[ |{\phi}|^2 \im {y} \,\Big|_u = |{\phi}|^2 \im {y} \,\Big|_{u_-}
+ \int_{u_-}^u \im V\, |{\phi}|^2 \]
and thus
\[ |{\phi}|^2 \im y \,\Big|_u \geq |{\phi}|^2 \im {y}  \,\Big|_{u_-} + \Big( \min_{[u_-, u_+]} \im V \Big)\:
\max_{[u_-,u_+]} |\phi|^2\: (u_+-u_-) \:. \]
Applying the Gr\"onwall estimate~\eqref{gronwall} gives the result.
\QED

\section{Semiclassical Estimates for a General Potential} \label{sec4}
\subsection{Estimates in the Case~$\re V < 0$} \label{sec41}
We now consider the Riccati equation~\eqref{riccati} on an interval~$I$.
We assume that the region~$I$ is semi-classical in the sense that the inequalities
\beq \label{CVbound}
\sup_I |V'| \leq \varepsilon\, \inf_I |V|^\frac{3}{2} \:,\quad
\sup_I |V''| \leq \varepsilon^2\, \inf_I |V|^2 \:,\quad
\sup_I |V'''| \leq \varepsilon^3\, \inf_I |V|^\frac{5}{2}
\eeq
hold, with a positive constant~$\varepsilon \ll 1$ to be specified later.

In this section, we derive estimates in the case~$\re V < 0$. As the approximate solution,
we choose the usual WKB wave function
\[ \tilde{\phi}(u) = V^{-\frac{1}{4}}\: \exp \Big( \int_{u_0}^u \sqrt{V} \Big) \:. \]
It is a solution of the Sturm-Liouville equation~\eqref{SLtilde} with
\begin{equation} \label{Vtilde}
\tilde{V} := V + \frac{5}{16}\: \frac{(V')^2}{V^2} - \frac{1}{4}\: \frac{V''}{V} \:.
\end{equation}
The corresponding solution of the Riccati equation~\eqref{RICtilde} becomes
\beq \label{yWKB}
\tilde{y} = \frac{\tilde{\phi}'}{\tilde{\phi}} = \sqrt{V} - \frac{V'}{4 V}\:.
\eeq
Moreover, we can compute the function~$\sigma$ from~\eqref{sigmarel},
\[ \sigma(u) = |\tilde{\phi}(u)|^2 = |V|^{-\frac{1}{2}} \,e^{2 \int_{u_0}^u \re \sqrt{V}} \:. \]

We begin with an estimate in the case~$\im V \geq 0$.
\begin{Lemma}  \label{lemmaWKBT2}
Assume that on the interval~$I:=[u_0, \umax]$, the potential~$V$ satisfies the
inequalities~\eqref{CVbound} with
\beq \label{epsbd2}
\varepsilon < \frac{1}{8}\:.
\eeq
Moreover, we assume that on~$I$,
\beq \label{ImRe2}
\im \sqrt{V} > \re \sqrt{V} \geq 0 \:.
\eeq
Then Theorem~\ref{thmT} applies with~$g \equiv 0$ and
\beq \label{Tcond}
\log T(u) \leq 64 \:\varepsilon^2\: \inf_I |V|^2 \int^u \frac{1}{|V|^\frac{3}{2}} \:.
\eeq
\end{Lemma}
\Proof The inequalities~\eqref{ImRe2} clearly imply that~$\im V \geq 0$.
Moreover, a straightforward calculation using~\eqref{Urel}, \eqref{tydef},
\eqref{yWKB} and~\eqref{Vtilde} shows that
\[ |U + \im^2 \sqrt{V}| \leq \frac{1}{2}\:
\frac{|V'|}{\sqrt{|V|}} + \frac{3}{8}\: \frac{|V'|^2}{|V|^2} + \frac{1}{4}\: \frac{|V''|}{|V|}
\leq 3 \varepsilon\, |V| \:, \]
where in the last step we used~\eqref{CVbound} and~\eqref{epsbd2}.
Combining this inequality with~\eqref{ImRe2} and~\eqref{epsbd2}, we conclude that
\[ U < -\frac{1}{4}\, |V| < 0 \:. \]
Hence Theorem~\ref{thmT} applies. Since~$\im V \geq 0$, we can satisfy the condition~\eqref{gdef}
by choosing~$g \equiv 0$.

A straightforward calculation and estimate (which we carried out with the help
of \textsf{Mathematica}) yields
\begin{gather}
|E_1+E_2| \leq 40\: \varepsilon^2\: \frac{\inf_I |V|^2}{|V|^\frac{3}{2}} \label{E12es} \\
|E_3| = \frac{|\im V|}{|U|}\: \frac{|\re(V-\tilde{V})|}{\sqrt{|U|} +\tilde{\beta}}
\leq 24\: \varepsilon^2\: \frac{\inf_I |V|^2}{|V|^\frac{3}{2}} \:, \label{E3es}
\end{gather}
giving the result.
\QED

The integral in~\eqref{Tcond} can be estimated efficiently if we assume that~$|V|$ satisfies
a weak version of concavity:
\begin{Lemma} \label{lemmaVabs}
Suppose that on the interval~$[u_0, u]$, the potential~$V$ satisfies the inequalities
\beq \label{absV}
|V(\tau)| \geq \frac{\tau-u_0}{u-u_0}\: |V(u)| + \frac{u-\tau}{u-u_0}\: |V(u_0)| \:.
\eeq
Then
\[ \int_{u_0}^u \frac{1}{|V|^\frac{3}{2}} \leq \frac{2\, (u - u_0)}{\sqrt{|V(u)|}\: |V(u_0)|}\:. \]
\end{Lemma}
\Proof Rewrite~\eqref{absV} as
\[ |V(\tau)| \geq |V(u)| + c\, (u - \tau) \qquad \text{with} \qquad
c := \frac{|V(u_0)| - |V(u)|}{u - u_0} \:. \]
Hence
\[ \int_{u_0}^u \frac{1}{|V|^\frac{3}{2}} \leq \int_{u_0}^\tau
\frac{d\tau}{(|V(u)| + c\, (u - \tau))^\frac{3}{2}} \:. \]
Computing and estimating the last integral gives the result.
\QED

The next lemma also applies in the case~$\im V <0$.
\begin{Lemma}  \label{lemmaWKBT3}
Assume that on the interval~$I:=[u_0, \umax]$, the potential~$V$ satisfies the
inequalities~\eqref{CVbound} with
\[ \varepsilon < \frac{1}{8}\:. \]
Moreover, we assume that for all~$u \in J=[u_0, u_1] \subset I$,
the inequalities~\eqref{absV} as well as the following inequalities hold:
\begin{align}
\im \sqrt{V} &> \re \sqrt{V} \geq 0 \label{ImRe21} \\
\sqrt{|V|} &\geq 200\, \varepsilon^2\, |J|\, \frac{\inf_I |V|^2}{|V(u_0)|} \label{Vl} \\
|J|\, |\im V| \:\sqrt{|V|} &\leq \frac{1}{30}\: |V(u_0)|\:. \label{Ju}
\end{align}
Then Theorem~\ref{thmT} applies on~$J$ if we choose~$T(u_0)=1$. Moreover,
\beq \label{Tes}
\log T \leq 100 \:\varepsilon^2\:\inf_I |V|^2 \int_{u_0}^u \frac{1}{|V|^\frac{3}{2}} \:.
\eeq
\end{Lemma}
\Proof The only difference to the proof of Lemma~\ref{lemmaWKBT2} is
that in order to satisfy~\eqref{gdef} we need to choose~$g$ positive.
Then the error term~$E_4$ is non-trivial. It is estimated by
\[ |E_4| \leq \frac{2\, |\im V|}{|V|^\frac{1}{2}} \:g \:. \]
In order to make this error term of about the same size as~\eqref{E12es}
and~\eqref{E3es}, we choose
\beq \label{grel}
g = 18\: \varepsilon^2\: \frac{\inf_I |V|^2}{|V|\: |\im V|} \:.
\eeq
Then the function~$T$ is bounded by~\eqref{Tes}.

Let us verify that the inequality~\eqref{gdef} is satisfied. Applying Lemma~\eqref{lemmaVabs},
we obtain
\begin{align*}
\log T \leq 200 \:\varepsilon^2\:\inf_I |V|^2 
\frac{|J|}{\sqrt{|V(u)|}\: |V(u_0)|} \:.
\end{align*}
Using~\eqref{Vl}, we see that the last expression is bounded by one. Hence, using the
mean value theorem,
\[ T-1 \leq e \: 200 \:\varepsilon^2\:\inf_I |V|^2 
\frac{|J|}{\sqrt{|V(u)|}\: |V(u_0)|} \:. \]
Comparing with~\eqref{grel} and using~\eqref{Ju}, we conclude that~\eqref{gdef} holds.
\QED

\subsection{Estimates in the Case~$\re V > 0$} \label{sec42}
We proceed with estimates in the case~$\re V > 0$.
We again assume that the inequalities~\eqref{CVbound} hold on an interval~$I$
for a suitable parameter~$\varepsilon>0$.
For the approximate solution~$\tilde{\phi}$, we now take the ansatz
\beq \label{tpansatz}
\tilde{\phi}(u) = V(u)^{-\frac{1}{4}}\: \exp \Big( \int_{0}^u \sqrt{V} + f \Big)
\eeq
with a so-called {\em{driving function}}~$f$ given by
\beq \label{fdef}
f := -\frac{s \varepsilon}{2}\:(1+i)\, \re \sqrt{V}
\eeq
and~$s \in \{-1, 1\}$. The function~$\tilde{\phi}$ is a solution of the Sturm-Liouville
equation~\eqref{SLtilde} with
\begin{equation} \label{Vtilde2}
\tilde{V} := (\sqrt{V}+f)^2 + \frac{5}{16}\: \frac{(V')^2}{V^2} - \frac{1}{4}\: \frac{V''}{V} 
-\frac{f}{2}\: \frac{V'}{V} + f' \:.
\end{equation}
The corresponding solution of the Riccati equation~\eqref{tydef} becomes
\beq \label{yWKB2}
\tilde{y} = \frac{\tilde{\phi}'}{\tilde{\phi}} = \sqrt{V} - \frac{V'}{4 V} + f\:.
\eeq
Again, we can compute the function~$\sigma$ from~\eqref{sigmarel} to obtain
\beq
\sigma = \frac{1}{\sqrt{|V|}} \: \exp \Big( 2 \,\re \int_{u_0}^u \sqrt{V} + f  \Big) 
\overset{\eqref{fdef}}{=} \frac{1}{\sqrt{|V|}} \: \exp \Big( (2 - s \varepsilon) \int_{u_0}^u \re \sqrt{V} \Big)  \:.
\label{sigmarel2}
\eeq

We want to apply the $\kappa$-method as introduced in Section~\ref{seckappa}.
We always choose~$\kappa(u_0)$ in agreement with~\eqref{tbeq}.
Again, in the symbols~$\pm$ and~$\mp$ the upper and lower case
refer to the cases~{\bf{(A)}} case~{\bf{(B)}}, respectively.
\begin{Lemma} \label{lemmaWKB1} 
Assume that on the interval~$I:=[u_0, \umax]$, the potential~$V$
satisfies the inequalities~\eqref{CVbound} with
\beq \label{eps2}
\varepsilon < \frac{1}{8} \:.
\eeq
Moreover, assume that
\beq
\big| \im \sqrt{V} \big| \leq \frac{1}{8}\:\re \sqrt{V} \:. \label{sqrtV}
\eeq
For a given parameter~$s \in \{1, -1\}$, we choose the approximate solution~$\tilde{\phi}$
of the form~\eqref{tpansatz} and~\eqref{fdef}. Then for all~$u \in I$, the following inequalities hold:
\begin{gather}
\frac{\varepsilon}{2}\:|V| \leq s \,\Big( U + \im^2 \sqrt{V} \Big) \leq 2 \varepsilon\:|V| \label{signU} \\
\varepsilon\, |V|^\frac{3}{2} \leq s \,\Big( {\mathfrak{D}} - (\kappa \pm R) \, \im V \Big) \leq 
3 \varepsilon\, |V|^\frac{3}{2} \label{Des} \\
\frac{1}{\sigma} \int_{u_0}^u \sigma\: |\im (V - \tilde{V})|
\leq 3 \varepsilon \:\sqrt{|V|}\:. \label{sigmaes}
\end{gather}
\end{Lemma}
\Proof Combining the identity
\[ |V| = \re^2 \sqrt{V} + \im^2 \sqrt{V} \]
with~\eqref{sqrtV}, we obtain
\begin{align}
\re^2 \sqrt{V} &\leq |V| \leq 2\: \re^2 \sqrt{V} \label{en1} \:.
\end{align}
Next, straightforward calculations using~\eqref{tpansatz}--\eqref{yWKB2} yield
\begin{align}
\im (V - \tilde{V}) &= s \varepsilon \re \sqrt{V}\, \left(  \re \sqrt{V}
+ \im \sqrt{V} \right) + E_1 \label{l2} \\
{\mathfrak{D}}\;\: &\!\!\!\overset{\eqref{Ddef3}}{=} s \varepsilon\, \re \sqrt{V}\,
\left\{ 2 \re^2 \sqrt{V} - \re \sqrt{V}\: \im \sqrt{V} + \im^2\, \sqrt{V} \right\} \notag \\
& \qquad + (\kappa \pm R)\, \im V + E_2 \label{l3} \\
U &\!\!\overset{\eqref{Udef}}{=} s \varepsilon \re^2 \sqrt{V} - \im^2 \sqrt{V} 
+ E_3 \:, \label{l5}
\end{align}
where the error terms~$E_1$, $E_2$ and~$E_3$ are estimated by
\begin{align}
|E_1| &\leq \frac{\varepsilon^2}{2}\: |V| + \varepsilon\: \frac{|V'|}{\sqrt{|V|}} +
\frac{5}{16}\: \frac{|V'|^2}{|V|^2} + \frac{1}{4}\: \frac{|V''|}{|V|}
\overset{\eqref{CVbound}}{\leq} 5 \varepsilon^2\, |V|  \label{l4} \\
|E_2| &\leq 9\: \frac{|V'|^2}{|V|^\frac{3}{2}} 
+\frac{9}{2}\: \frac{|V''|}{\sqrt{|V|}} 
+\frac{|V'''|}{|V|}
+ \frac{51}{8}\: \frac{|V'|^3}{|V|^3}
+ \frac{21}{4}\: \frac{|V'\, V''|}{|V|^2} \notag \\
&\quad+ (12 \,\varepsilon + 3 \varepsilon^2)\: |V'|
+ (12\, \varepsilon^2 + 2 \varepsilon^3)\: |V|^\frac{3}{2}
+ \frac{9 \varepsilon}{2}\: \frac{|V'|^2}{|V|^\frac{3}{2}}
+ \frac{9 \varepsilon}{4}\: \frac{|V''|}{\sqrt{|V|}} \notag \\
&\!\!\overset{\eqref{CVbound}}{\leq} 40\: \varepsilon^2\, |V|^{\frac{3}{2}} + 25\: \varepsilon^3\, |V|^{\frac{3}{2}}
\overset{\eqref{eps2}}{\leq} 50\: \varepsilon^2\, |V|^{\frac{3}{2}} \label{l52} \\
|E_3| &\leq \big|\im \sqrt{V} \big|\: \frac{|V'|}{|V|}
+ \frac{7}{4}\: \frac{|V'|^2}{|V|^2} + \frac{|V''|}{|V|} + \frac{3 \varepsilon}{4}\: \frac{|V'|}{\sqrt{|V|}}
+ \frac{\varepsilon^2}{4}\: |V| \notag \\
&\!\!\!\!\!\!\!\!\!\overset{\eqref{sqrtV}, \eqref{CVbound}}{\leq}
\frac{\varepsilon}{15}\:|V| + 2\: \varepsilon^2\: |V|
\overset{\eqref{eps2}}{\leq} \frac{\varepsilon}{3}\:|V| \:. \label{l6}
\end{align}

The estimate~\eqref{signU} follows immediately from~\eqref{l5} and~\eqref{l6}
combined with~\eqref{en1} and~\eqref{eps2}.

In order to prove~\eqref{Des}, we estimate the curly brackets in~\eqref{l3} from above and below
using the Schwarz inequality,
\begin{align*}
2 \re^2 \sqrt{V} & - \re \sqrt{V}\: \im \sqrt{V} + \im^2 \sqrt{V} \leq
\frac{5}{2} \left( \re^2 \sqrt{V} + \im^2 \sqrt{V} \right) = \frac{5}{2}\:|V| \\
2 \re^2 \sqrt{V} & - \re \sqrt{V}\: \im \sqrt{V} + \im^2 \sqrt{V} \geq
\frac{3}{2}\: \re^2 \sqrt{V} + \frac{1}{2}\: \im^2 \sqrt{V}
\overset{\eqref{sqrtV}}{\geq} \frac{5}{4}\: |V|\:.
\end{align*}
Using~\eqref{eps2} in~\eqref{l6}, we can compensate the error term~$E_2$ in~\eqref{l3} to obtain~\eqref{Des}.

It remains to prove~\eqref{sigmaes}: We first apply~\eqref{l2} and~\eqref{l4} to obtain
\begin{align*}
\frac{1}{\sigma} \int_{u_0}^u \sigma\: |\im (V - \tilde{V})|
&\leq \frac{\varepsilon}{\sigma} \int_{u_0}^u \sigma\: \sqrt{|V|}
\left( \re \sqrt{V} + \im \sqrt{V} \right) \leq \frac{2 \varepsilon}{\sigma}
\int_{u_0}^u \sigma\: \sqrt{|V|} \:\re \sqrt{V} \:.
\end{align*}
Using~\eqref{sigmarel2}, we obtain
\begin{align*}
\frac{1}{\sigma} \int_{u_0}^u &\sigma\: |\im (V - \tilde{V})|
\leq \frac{2 \varepsilon}{\sigma} \int_{u_0}^u e^{ (2 - s \varepsilon) \int_{u_0}^v \re \sqrt{V} } \:\re \sqrt{V(v)}\:dv \\
&\leq \frac{2 \varepsilon}{\sigma \,(2 - s \varepsilon)}
\int_{u_0}^u \frac{d}{dv}
e^{ (2 - s \varepsilon) \,\re \int_{u_0}^v \sqrt{V} }\:dv 
= \frac{2 \varepsilon}{\sigma\, (2 - s \varepsilon)}
\left( e^{ (2 - s \varepsilon) \int_{u_0}^u \re \sqrt{V} } - 1 \right) \\
&\!\!\!\overset{\eqref{sigmarel2}}{=}  \frac{2 \varepsilon}{2 - s \varepsilon} \:\sqrt{|V|} \left( 1 -
e^{ -(2 - s \varepsilon) \int_{u_0}^u \re \sqrt{V} } \right) .
\end{align*}
Applying~\eqref{eps2} and~\eqref{sqrtV} gives~\eqref{sigmaes}.
\QED

So far, we did not specify the function~$\kappa$. If~$\kappa$ is chosen according to~\eqref{tbeq2},
then one can apply Theorem~\ref{thm1} in both case~{\bf{(A)}} or~{\bf{(B)}}, provided that the
determinator has the correct sign.
We now explore the possibilities for applying Proposition~\ref{prpnewkappa}.
\begin{Lemma} \label{lemmag}
Suppose that the function~$g$ in~\eqref{kset} is chosen as
\[ g = \nu\, \sqrt{|V|}\, \sigma \:, \]
where the positive parameters~$\varepsilon$ and~$\nu$ satisfy the following conditions,
\begin{gather}
100 \,\varepsilon^2 < \nu^2 < \varepsilon < \frac{1}{100} \label{epsdel} \\
|\im \sqrt{V}| \leq \frac{\nu}{10}\: \re \sqrt{V} \:. \label{imVupper}
\end{gather}
Then the function~$g$ is monotone increasing. Choosing again the ansatz~\eqref{tpansatz} with the driving
function~\eqref{fdef} and~$s=1$, the determinator is positive. Moreover,
\begin{align}
\big| \tilde{\beta} + \kappa \big| &\leq \frac{3}{2}\: \nu\, \sqrt{|V|} \label{bkupper} \\
\big| (\kappa - R)\, \im V \big| &\leq \frac{3}{5}\:\varepsilon\: |V|^\frac{3}{2} \:. \label{imV2}
\end{align}
\end{Lemma}
\Proof We first note that the assumptions~\eqref{epsdel} and~\eqref{imVupper}
imply that~\eqref{eps2} and~\eqref{sqrtV} are satisfied, so that we may use Lemma~\ref{lemmaWKB1}.
According to~\eqref{sigmarel2},
\[ g = \nu\, \exp \Big( (2 - s \varepsilon) \int_{u_0}^u \re \sqrt{V} \Big) \:, \]
which is indeed increasing in view of~\eqref{imVupper}.
Next, according to~\eqref{kset},
\[ \kappa(u) = \nu\, \sqrt{|V|} + \frac{1}{\sigma} \int_{u_0}^u \sigma \im(V-\tilde{V}) \:. \]
In view of~\eqref{epsdel}, we know that~$\nu > 10 \varepsilon$. Also using the estimate~\eqref{sigmaes},
one finds that
\beq \label{kupper}
\sqrt{|V|} \:(\nu - 3 \varepsilon) \leq \kappa(u) \leq \sqrt{|V|} \:(\nu + 3 \varepsilon) \:.
\eeq
Moreover, using~\eqref{yWKB2} and~\eqref{CVbound},
\begin{align*}
\tilde{\beta} &= \im \sqrt{V} - \im \frac{V'}{4 V} - \frac{s \varepsilon}{2}\: \re \sqrt{V} \\
|\tilde{\beta}| &\leq \frac{1}{40} \: \sqrt{|V|} \left( 30 \varepsilon + 4 \nu \right)
\end{align*}
and thus, using~\eqref{epsdel} and~\eqref{kupper},
\begin{align}
\big| \tilde{\beta} + \kappa - \nu \sqrt{|V|} \big| &\leq \frac{\nu}{2}\: \sqrt{|V|} \label{bk2} \\
\big| \tilde{\beta} + \kappa \big| &\geq \frac{\nu}{2}\: \sqrt{|V|} \:. \label{bklower}
\end{align}
Moreover, \eqref{bk2} yields~\eqref{bkupper}.

We next apply Lemma~\ref{lemmaRpmb}. Combining~\eqref{kupper} and~\eqref{bklower}
with~\eqref{Vtilde2} and~\eqref{CVbound}, we can use~\eqref{epsdel} to obtain
\[ |\kappa - R| \leq 3 \varepsilon\: \frac{\sqrt{|V|}}{\nu} \qquad \text{and} \qquad
|\im V| \leq \frac{\nu}{5}\: |V| \:. \]
This proves~\eqref{imV2}. Using this inequality in~\eqref{Des} concludes the proof.
\QED

\section{Semiclassical Estimates for the Angular Teukolsky Equation} \label{sec5}
\subsection{Estimates in the Case~$\re V < 0$}
\begin{figure}
\begin{picture}(0,0)%
\includegraphics{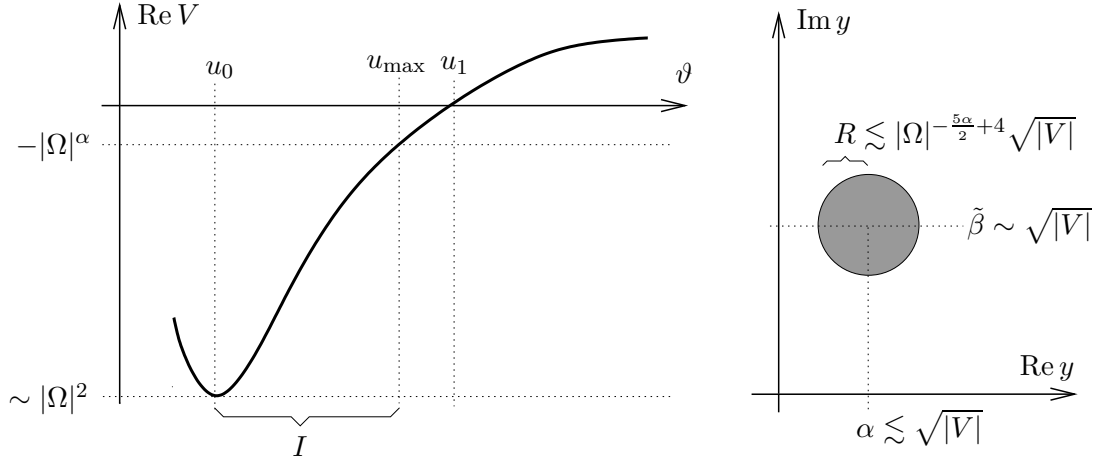}%
\end{picture}%
\setlength{\unitlength}{2486sp}%
\begingroup\makeatletter\ifx\SetFigFont\undefined%
\gdef\SetFigFont#1#2#3#4#5{%
  \reset@font\fontsize{#1}{#2pt}%
  \fontfamily{#3}\fontseries{#4}\fontshape{#5}%
  \selectfont}%
\fi\endgroup%
\begin{picture}(10696,4646)(320,-8060)
\put(1621,-3661){\makebox(0,0)[lb]{\smash{{\SetFigFont{11}{13.2}{\familydefault}{\mddefault}{\updefault}${\text{Re}}\,V$}}}}
\put(2296,-4156){\makebox(0,0)[lb]{\smash{{\SetFigFont{11}{13.2}{\familydefault}{\mddefault}{\updefault}$u_0$}}}}
\put(3166,-7973){\makebox(0,0)[lb]{\smash{{\SetFigFont{11}{13.2}{\rmdefault}{\mddefault}{\updefault}$I$}}}}
\put(406,-4951){\makebox(0,0)[lb]{\smash{{\SetFigFont{11}{13.2}{\familydefault}{\mddefault}{\updefault}$-|\Omega|^\alpha$}}}}
\put(335,-7482){\makebox(0,0)[lb]{\smash{{\SetFigFont{11}{13.2}{\rmdefault}{\mddefault}{\updefault}$\sim |\Omega|^2$}}}}
\put(6984,-4246){\makebox(0,0)[lb]{\smash{{\SetFigFont{11}{13.2}{\familydefault}{\mddefault}{\updefault}$\vartheta$}}}}
\put(3901,-4096){\makebox(0,0)[lb]{\smash{{\SetFigFont{11}{13.2}{\familydefault}{\mddefault}{\updefault}$u_\text{max}$}}}}
\put(4636,-4111){\makebox(0,0)[lb]{\smash{{\SetFigFont{11}{13.2}{\familydefault}{\mddefault}{\updefault}$u_1$}}}}
\put(8191,-3706){\makebox(0,0)[lb]{\smash{{\SetFigFont{11}{13.2}{\familydefault}{\mddefault}{\updefault}${\text{Im}}\,y$}}}}
\put(9890,-5759){\makebox(0,0)[lb]{\smash{{\SetFigFont{11}{13.2}{\familydefault}{\mddefault}{\updefault}$\tilde{\beta} \sim \sqrt{|V|}$}}}}
\put(8560,-4848){\makebox(0,0)[lb]{\smash{{\SetFigFont{11}{13.2}{\familydefault}{\mddefault}{\updefault}$R \lesssim |\Omega|^{-\frac{5 \alpha}{2}+4} \sqrt{|V|}$}}}}
\put(8787,-7779){\makebox(0,0)[lb]{\smash{{\SetFigFont{11}{13.2}{\familydefault}{\mddefault}{\updefault}$\alpha \lesssim \sqrt{|V|}$}}}}
\put(10419,-7167){\makebox(0,0)[lb]{\smash{{\SetFigFont{11}{13.2}{\familydefault}{\mddefault}{\updefault}${\text{Re}}\,y$}}}}
\end{picture}%
\caption{WKB estimate in the case~$\re V < 0$.}
\label{figWKB1}
\end{figure}

We now apply the estimates of Section~\ref{sec41} to the
angular Teukolsky equation. We choose~$u_0$ and~$u_1$ as the minimum
and the zero of the real part of the potential, respectively,
\[ \re V'(u_0) = 0 \qquad \text{and} \qquad \re V(u_1) = 0 \]
(see the left of Figure~\ref{figWKB1}). 
In order to simplify the notation in our estimates we use the notation
\[ f \lesssim |\Omega|^\beta \qquad \text{for the inequality} \qquad
|f| \leq \text{c}\, |\Omega|^\beta \]
with a constant~$c$ which is independent of the parameters~$\Omega$ and~$\mu$
under consideration. Likewise, we use the symbol
\[ f \eqsim |\Omega|^\beta \qquad \text{for} \qquad
\frac{1}{c}\:|\Omega|^\beta \leq |f| \leq c\: |\Omega|^\beta \:. \]
We choose~$\umax$ such that
\[ \re V(u_\text{max}) = -|\Omega|^\alpha \qquad \text{with~$1 < \alpha < 2$}\:. \]

We now prove the invariant disk estimate illustrated on the right of Figure~\ref{figWKB1}.
\begin{Prp} \label{prpWKBm}
For any~$\alpha$ in the range
\[ \frac{8}{5} < \alpha \leq 2 \]
and sufficiently large~$C$,
we consider the invariant region estimate of Theorem~\ref{thm2}
on the interval~$I=[u_0, u_\text{max}]$ with the initial condition~$T(u_0)=1$,
taking the WKB solution~\eqref{yWKB} as our approximate solution.
Moreover, we consider~$\Omega$ of the form~\eqref{Obounds} such that
\[ \im V|_I \geq 0\:. \]
Then the invariant region estimate applies on~$I$, and the function~$T$ is bounded by
\[ \log T(u) \lesssim |\Omega|^{-\frac{5 \alpha}{2} + 4}\:. \]
\end{Prp}
\Proof We want to apply Lemma~\ref{lemmaWKBT2}. We choose
\[ \varepsilon = |\Omega|^{2-\frac{3 \alpha}{2}}\:. \]
The estimates
\begin{align*}
\sup_I |V'| &\lesssim |\Omega|^2 = \varepsilon\, |\Omega|^\frac{3 \alpha}{2} 
\lesssim \varepsilon \inf_I |V|^\frac{3}{2} \\
\sup_I |V''| &\lesssim |\Omega|^2 \lesssim  \varepsilon^2\, |\Omega|^{2 \alpha} \simeq \varepsilon^2 \inf_I |V|^2  \\
\sup_I |V'''| &\lesssim |\Omega|^2 \lesssim \varepsilon^3 \,|\Omega|^\frac{5 \alpha}{2}
\simeq \varepsilon^3 \inf_I |V|^\frac{5}{2} 
\end{align*}
show that for that for large~$|\Omega|$, the WKB conditions~\eqref{CVbound} hold.

In order to verify~\eqref{ImRe2}, we note that the inequalities~$\re V \lesssim |\Omega|^\alpha$
and~$0 \leq \im V \lesssim |\Omega|$ imply that the argument of~$V$ lies in the
interval~$[150^\circ, 180^\circ)$. Choosing the sign convention for the square root
such that~$\arg \sqrt{V} \in [75^\circ, 90^\circ)$, proving~\eqref{ImRe2}.

We finally estimate~\eqref{Tcond} by
\[ \varepsilon^2\, \inf_I |V|^2 \int_{u_0}^u \frac{1}{|V|^\frac{3}{2}} \leq
|I|\, \varepsilon^2\, \inf_I \sqrt{|V|} \lesssim
|\Omega|^{4 - 3 \alpha}\: |\Omega|^\frac{\alpha}{2}
= |\Omega|^{-\frac{5 \alpha}{2}+4} \:, \]
concluding the proof.
\QED

\subsection{Estimates in the Case~$\re V > 0$}

\begin{figure}
\begin{picture}(0,0)%
\includegraphics{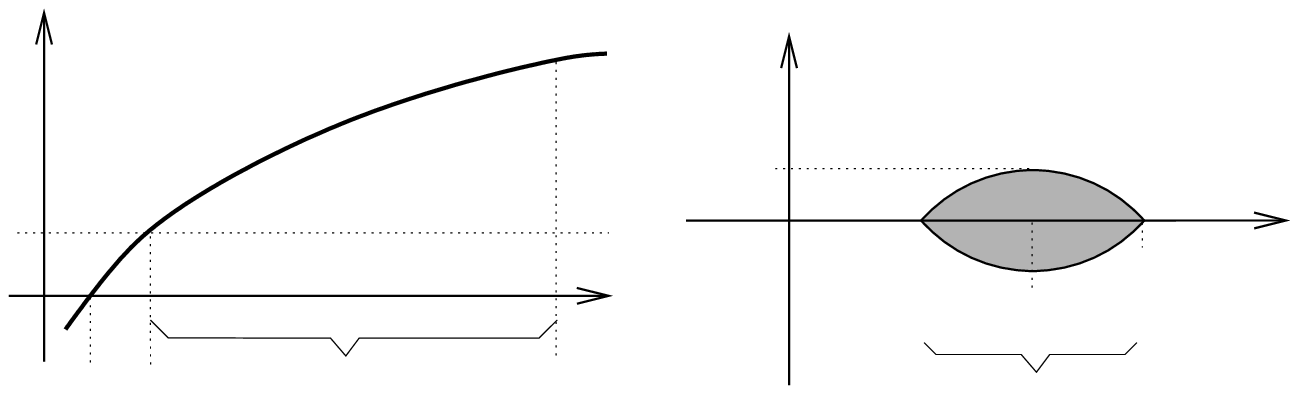}%
\end{picture}%
\setlength{\unitlength}{2486sp}%
\begingroup\makeatletter\ifx\SetFigFont\undefined%
\gdef\SetFigFont#1#2#3#4#5{%
  \reset@font\fontsize{#1}{#2pt}%
  \fontfamily{#3}\fontseries{#4}\fontshape{#5}%
  \selectfont}%
\fi\endgroup%
\begin{picture}(10606,3275)(800,-7927)
\put(2101,-7689){\makebox(0,0)[lb]{\smash{{\SetFigFont{12}{14.4}{\familydefault}{\mddefault}{\updefault}$u_0$}}}}
\put(2559,-7704){\makebox(0,0)[lb]{\smash{{\SetFigFont{11}{13.2}{\familydefault}{\mddefault}{\updefault}$u_\text{min}$}}}}
\put(6068,-6833){\makebox(0,0)[lb]{\smash{{\SetFigFont{11}{13.2}{\familydefault}{\mddefault}{\updefault}$\vartheta$}}}}
\put(1973,-4883){\makebox(0,0)[lb]{\smash{{\SetFigFont{11}{13.2}{\familydefault}{\mddefault}{\updefault}${\text{Re}}\,V$}}}}
\put(5630,-7643){\makebox(0,0)[lb]{\smash{{\SetFigFont{11}{13.2}{\familydefault}{\mddefault}{\updefault}$u_1$}}}}
\put(9893,-6826){\makebox(0,0)[lb]{\smash{{\SetFigFont{11}{13.2}{\familydefault}{\mddefault}{\updefault}$\alpha + \sqrt{U}$}}}}
\put(9331,-7186){\makebox(0,0)[lb]{\smash{{\SetFigFont{11}{13.2}{\familydefault}{\mddefault}{\updefault}$\alpha \eqsim \text{Re}\, V$}}}}
\put(11064,-6165){\makebox(0,0)[lb]{\smash{{\SetFigFont{11}{13.2}{\familydefault}{\mddefault}{\updefault}${\text{Re}}\,y$}}}}
\put(7920,-5581){\makebox(0,0)[lb]{\smash{{\SetFigFont{11}{13.2}{\familydefault}{\mddefault}{\updefault}$\beta+R \eqsim |\Omega|^{2-\frac{3 \alpha}{2}}\, \sqrt{|V|}$}}}}
\put(7711,-5131){\makebox(0,0)[lb]{\smash{{\SetFigFont{11}{13.2}{\familydefault}{\mddefault}{\updefault}${\text{Im}}\,y$}}}}
\put(6629,-6062){\makebox(0,0)[lb]{\smash{{\SetFigFont{11}{13.2}{\familydefault}{\mddefault}{\updefault}$\beta+R$}}}}
\put(9308,-7831){\makebox(0,0)[lb]{\smash{{\SetFigFont{11}{13.2}{\familydefault}{\mddefault}{\updefault}$\eqsim |\Omega|^{1-\frac{3 \alpha}{4}}\, \sqrt{|V|}$}}}}
\put(4115,-7711){\makebox(0,0)[lb]{\smash{{\SetFigFont{11}{13.2}{\rmdefault}{\mddefault}{\updefault}$I$}}}}
\put(815,-6564){\makebox(0,0)[lb]{\smash{{\SetFigFont{11}{13.2}{\familydefault}{\mddefault}{\updefault}$C\, |\Omega|^\alpha$}}}}
\end{picture}%
\caption{WKB estimate in the case~$\re V > 0$.}
\label{figWKB2}
\end{figure}

In order to apply Lemma~\ref{lemmaWKB1}, we consider~$u_0$ such that
\beq \label{u0def}
\re V(u_0) = 0 \qquad \text{and} \qquad \re V'(u_0) \eqsim |\Omega|^2\:.
\eeq
We choose~$u_\text{min} > u_0$ such that
\[ \re V(u_\text{min}) = C \,|\Omega|^\alpha \]
with
\beq \label{alpharange2}
\frac{4}{3} \leq \alpha < 2
\eeq
and a constant~$C$ to be chosen independent of~$\Omega$
(see the left of Figure~\ref{figWKB2}).
Moreover, we assume that
\beq 
\re V|_I \geq \frac{C}{2}\: |\Omega|^\alpha\:. \label{Vgtr}
\eeq
We next apply the invariant region estimates of Proposition~\ref{prpnewkappa},
relying on the estimates of Lemmas~\ref{lemmaWKB1} and~\ref{lemmag}.
We introduce the set~${\mathcal{I}}$ as the intersection of the upper half plane with the
circle with center~$m=\alpha + i \beta$ and radius~$R$,
\[ {\mathcal{I}} = \left\{ z \in \C \:|\: |z-m| \leq R \text{ and } \re z \geq 0 \right\} \]
(where again~$\alpha=\re \tilde{y}$ and~$\tilde{y}$ as in~\eqref{yWKB2} and~\eqref{fdef}).
Moreover, we let~$\overline{\mathcal{I}}$ be the complex conjugate of the set~${\mathcal{I}}$.

\begin{Prp} \label{prpWKB2}
We choose the interval~$I=[u_\text{min}, u_1]$ according to~\eqref{u0def}--\eqref{Vgtr}.
Assume that that~$V$ satisfies on~$I$ the conditions~\eqref{CVbound} and~\eqref{sqrtV}.
Then the region~${\mathcal{I}} \cup \overline{\mathcal{I}}$ is invariant under the Riccati flow.
Moreover,
\beq \label{WKBlens}
R + \beta \eqsim |\Omega|^{2-\frac{3 \alpha}{2}}\, \sqrt{|V|}\:,\qquad
|R^2 - \beta^2| \lesssim |\Omega|^{2-\frac{3 \alpha}{2}}\: |V|\:.
\eeq
\end{Prp} \noindent
Before giving the proof, we note that in the case~$U<0$, the sets~${\mathcal{I}}$
and~$\overline{\mathcal{I}}$ do not intersect, so that the invariant region are two disjoint disks.
In the case~$U>0$, the two disks form a connected set.
In the case~$\beta<0$, we obtain a {\em{lens-shaped invariant region}}, as
as illustrated on the right of Figure~\ref{figWKB2}.
\Proof[Proof of Proposition~\ref{prpWKB2}]
Similar as in the proof of Proposition~\ref{prpWKBm}, a Taylor expansion of the potential
around~$u_0$ yields that
\[ u_\text{min} - u_0 \eqsim |\Omega|^{\alpha-2}\:. \]
We want to choose~$\varepsilon$ as small as possible, but in agreement with~\eqref{CVbound}.
This leads us to make the ansatz
\[ \varepsilon = \delta\, |\Omega|^{2 - \frac{3 \alpha}{2}} \]
with~$0 < \delta \ll 1$ independent of~$|\Omega|$.
By choosing~$\delta$ sufficiently small and~$C$ sufficiently large,
we can arrange that the inequalities~\eqref{CVbound}, \eqref{eps2}
as well as the last inequality in~\eqref{epsdel} hold.
Next we choose~$\nu$ in agreement with~\eqref{epsdel}, but as small as possible,
\[ \nu = 20\, \delta\, |\Omega|^{2-\frac{3 \alpha}{2}} \:. \]

Let us verify that Lemmas~\ref{lemmaWKB1} and~\ref{lemmag} apply.
As just explained, \eqref{eps2} holds for sufficiently small~$\delta$.
According to~\eqref{Obounds}, we know that
\[ |\Omega| \gtrsim |\im V| = 2 \re \sqrt{V} \,\im \sqrt{V} \]
and thus in view of~\eqref{Vgtr},
\beq \label{scal2}
\re \sqrt{V} \gtrsim C^\frac{1}{2}\, |\Omega|^{\frac{\alpha}{2}} \qquad \text{and} \qquad
\big| \im \sqrt{V} \big| \lesssim C^{-\frac{1}{2}} \,|\Omega|^{1-\frac{\alpha}{2}}\:.
\eeq
Thus, possibly after increasing~$C$, the inequality~\eqref{sqrtV} is satisfied.
Hence Lemma~\ref{lemmaWKB1} applies.
The inequalities~\eqref{epsdel} again hold for sufficiently small~$\delta$.
Using~\eqref{scal2}, we see that
\[ \left| \frac{\im \sqrt{V}}{\re \sqrt{V}} \right| \lesssim C^{-1}\, |\Omega|^{1-\alpha}
= \frac{\nu}{10} \: \frac{|\Omega|^{-1+\frac{\alpha}{2}}}{2 \delta \,C} \:, \]
and in view of~\eqref{alpharange2}, the last factor can be made arbitrarily small
by further increasing~$C$ if necessary.
Hence~\eqref{imVupper} holds, and Lemma~\ref{lemmag} applies.

We begin with the case when~$y$ lies in the upper half plane (the general case will be treated below).
Choosing~$s=1$, we can apply Proposition~\ref{prpnewkappa}
to obtain the invariant region estimate~\eqref{ninvdisk}.
The first inequality in~\eqref{WKBlens} follows from the first equation in~\eqref{ninvdisk}
and~\eqref{bkupper}.
Similarly, the second inequality in~\eqref{WKBlens} follows from the second equation in~\eqref{ninvdisk}
and~\eqref{signU}, noting that according to~\eqref{scal2},
\[ \im^2 \sqrt{V} \lesssim |\Omega|^{2-\alpha} \lesssim |\Omega|^{2-2 \alpha} \:|V|
\lesssim \varepsilon |V| \:. \]

If~$y$ lies in the lower half plane, we take the complex conjugate of the Riccati equation
and again apply the above estimates. This simply amounts to flipping the sign of~$\beta$ in all formulas.
If~$y(u)$ crosses the real line, we can perform the replacement~$\beta \rightarrow -\beta$,
which describes a reflection of the invariant circle at the real axis.
In this way, we can flip from estimates in the upper to estimates in the lower half plane
and vice versa, without violating our estimates. We conclude that~$y$ stays inside the lens-shaped region
obtained as the intersection of the two corresponding invariant circles.
\QED

\section{Parabolic Cylinder Estimates} \label{sec6}
Near the turning points of the real part of the potential, we approximate the potential
by a quadratic polynomial,
\beq \label{tVparabolic}
\tilde{V}(u) = {\mathfrak{p}} + \frac{{\mathfrak{q}}}{4}\: (u - {\mathfrak{r}})^2 \qquad \text{with} \qquad {\mathfrak{p}}, {\mathfrak{q}}, {\mathfrak{r}} \in \C\:.
\eeq
The corresponding differential equation~\eqref{SLtilde} can be solved explicitly in terms
of the parabolic cylinder function, as we now recall. The parabolic cylinder function,
which we denote by~$U_a(z)$, is a solution of the differential equation
\[ U_a''(z) = \Big( \frac{z^2}{4}\: +a \Big) \, U_a(z)\:. \]
Setting
\beq \label{azdef}
\tilde{\phi}(u) = U_a(z) \qquad \text{with} \qquad
a = \frac{{\mathfrak{p}}}{\sqrt{{\mathfrak{q}}}} \:,\quad z = {\mathfrak{q}}^{\frac{1}{4}}\, (u-{\mathfrak{r}}) \:,
\eeq
a short calculation shows that~$\tilde{\phi}$ indeed satisfies~\eqref{SLtilde}.
We set
\[ b=-4 \,\Big(a-\frac{1}{2} \Big) \:. \]

\subsection{Estimates of Parabolic Cylinder Functions}
In preparation for getting invariant region estimates, we need to
get good control of the parabolic cylinder function~$U_a(z)$.
To this end, in this section we elaborate on the general results in~\cite{special}
and bring them into a form which is most convenient for our applications.

\begin{Lemma} 
There is a constant~$c>0$ such that for all parameters~$z$, $b$ in the range
\[ |z|^2 > c \qquad \text{and} \qquad |z|^2 > 4 |b| \:, \]
the parabolic cylinder function is well-approximated by the WKB solution.
\end{Lemma}
\Proof We want to apply~\cite[Theorem~3.3]{special} in the case~$t_0 = t_+$
(with~$t_+$ as defined in~\cite[eqn~(3.10)]{special}). 
Using~\cite[eqns~(3.14) and~(3.17)]{special}, we find
\[ |8 d| \geq \left| z + \sqrt{z^2 - b} \right|^2 \geq 
\left| 2 \sqrt{|b|} + \sqrt{3 |b|} \right|^2 = (2 + \sqrt{3})^2 \:|b| > 8 \,|b| \:. \]
Hence the parameter~$\rho$ defined in~\cite[eqn~(3.17)]{special}) is smaller than~$1/8$,
making it possible to choose~$\kappa=1/4$ (see~\cite[Lemma~3.2]{special}).
Applying~\cite[Theorem~3.3]{special} gives the result.
\QED

For the following estimates, we work with the Airy-WKB limit, giving us the asymptotic
solution~\cite[eqns~(3.36) and~(3.37)]{special}.
\begin{Lemma} Assume that 
\[ |z^2-b| \leq |b|^\frac{1}{3} \:, \qquad \arg b \in (88^\circ, 92^\circ)
\qquad \text{and} \qquad |b| > 100\:. \]
Then the estimate of~\cite[Theorem~3.9]{special} applies and~$|h(z)|^2<2$.
\end{Lemma}
\Proof According to~\cite[eqns~(3.10) and~(3.37)]{special}
\begin{align*}
4 t_0^2 - b &= 2 (z^2-b) \pm 2 z \sqrt{z^2-b} \\
|h(z)|^2 &= \frac{1}{4 \!\cdot\! 2^\frac{4}{3}} \: |b|^{-\frac{4}{3}} \,|4 t_0^2 - b|^2
\end{align*}
and thus
\begin{align}
|z|^2 &\leq |z^2-b| + |b| \leq |b| + |b|^{\frac{1}{3}} \notag \\
|z| &\leq |b|^\frac{1}{2} + |b|^{\frac{1}{6}} \notag \\
\big| z \pm \sqrt{z^2-b} \big| &\leq |z| + |b|^{\frac{1}{3}} \leq
|b|^\frac{1}{2} + 2\, |b|^{\frac{1}{6}} \notag \\
\big| z \pm \sqrt{z^2-b} \big|^2 &\leq |b| + 8\, |b|^{\frac{2}{3}} \notag \\
|4 t_0^2 - b|^2 &\leq 4 \,|z^2-b|\, \Big| z \pm \sqrt{z^2-b} \Big|^2 \notag \\
&\leq 4\, |b|^{\frac{4}{3}} \left( 1 + 8\, |b|^{-\frac{1}{3}} \right) \notag \\
|h(z)|^2 &\leq \frac{1}{2^\frac{4}{3}} \left( 1 + 8\, |b|^{-\frac{1}{3}} \right) < 2 \notag \\
\left| \frac{4 t_0^2}{b} - 1 \right|^2 &\leq
4\, |b|^{-\frac{2}{3}} \left( 1 + 8\, |b|^{-\frac{1}{3}} \right) < 0.6\:. \label{rest}
\end{align}
We now apply~\cite[Theorem~3.9]{special}, noting that~\eqref{rest} implies
the condition~\cite[eqn~(3.39)]{special}.
\QED

\begin{figure}
\begin{picture}(0,0)%
\includegraphics{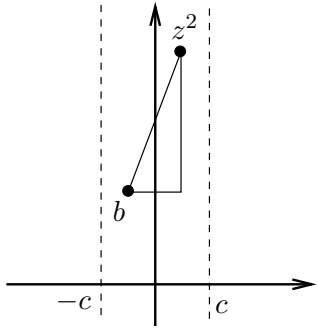}%
\end{picture}%
\setlength{\unitlength}{1989sp}%
\begingroup\makeatletter\ifx\SetFigFont\undefined%
\gdef\SetFigFont#1#2#3#4#5{%
  \reset@font\fontsize{#1}{#2pt}%
  \fontfamily{#3}\fontseries{#4}\fontshape{#5}%
  \selectfont}%
\fi\endgroup%
\begin{picture}(3936,4116)(363,-7294)
\put(2998,-7096){\makebox(0,0)[lb]{\smash{{\SetFigFont{11}{13.2}{\familydefault}{\mddefault}{\updefault}$c$}}}}
\put(1723,-5951){\makebox(0,0)[lb]{\smash{{\SetFigFont{11}{13.2}{\familydefault}{\mddefault}{\updefault}$b$}}}}
\put(2448,-3681){\makebox(0,0)[lb]{\smash{{\SetFigFont{11}{13.2}{\familydefault}{\mddefault}{\updefault}$z^2$}}}}
\put(1013,-7051){\makebox(0,0)[lb]{\smash{{\SetFigFont{11}{13.2}{\familydefault}{\mddefault}{\updefault}$-c$}}}}
\end{picture}%
\caption{Estimating the argument of~$z^2-b$.}
\label{figz2b}
\end{figure}

\begin{Lemma} \label{lemmaaway}
For any~$c>0$ there is a constant~${\mathscr{C}} >0$ such that the
following statement is valid. Assume that
\[ |z^2-b| > |b|^\frac{1}{3} \:,\qquad |\re z^2|, |\re b| < c \qquad \text{and} \qquad
\im z^2, \im b > {\mathscr{C}} \:. \]
Then the assumptions of~\cite[Theorem~3.9]{special} hold. Moreover,
the argument~$h^2$ of the Airy function in~\cite[eqn~(3.36)]{special}
avoids the branch cut (i.e.\ there is a constant~$\varepsilon(c, {\mathscr{C}})>0$
such that~\cite[eqn~(2.6)]{special} holds).
As a consequence, the Airy function has the WKB approximation given in~\cite[Theorem~2.2]{special}.
\end{Lemma}
\Proof By choosing~${\mathscr{C}}$ sufficiently large, we can arrange that the
arguments of~$z^2$ and~$b$ are arbitrarily close to~$90^\circ$.
Moreover, as shown in Figure~\ref{figz2b}, we have the inequality
\[ \cos \arg (z^2-b) \leq \frac{2c}{|z^2-b|} \leq 2 c\: |b|^{-\frac{1}{3}} \leq 2 c\, {\mathscr{C}}^{-\frac{1}{3}}\:, \]
showing that for sufficiently large~${\mathscr{C}}$, the argument of~$z^2-b$ is arbitrarily close
to~$\pm 90^\circ$.

We next consider the phase of~$t_0$ given by either~$t_+$ or~$t_-$,
\beq \label{t0def}
t_0 = \frac{1}{2} \left(z \pm \sqrt{z^2 -b} \right) .
\eeq
We need to consider both signs in order to take into account both branches of the square root.
Since the arguments of both~$z^2$ and~$z^2-b$ are arbitrarily close to~$90^\circ$, we
know that the arguments of~$z$ and~$\sqrt{z^2-b}$ are both arbitrarily close to~$45^\circ \!\!\!\! \mod 180^\circ$.
Hence choosing the sign in~\eqref{t0def} such that the real parts of~$z$ and~$\pm \sqrt{z^2-b}$
have the signs, it follows immediately that the argument of~$t_0$ is also
arbitrarily close to~$45^\circ \!\!\!\! \mod 180^\circ$. The identity
\[ t_+ \,t_- = \frac{b}{4} \]
yields that for sufficiently large~${\mathcal{C}}$, the argument of the other branch
is also arbitrarily close to~$45^\circ \!\!\!\! \mod 180^\circ$.

As a consequence, the conditions~\cite[eqns~(3.38) and~(3.39)]{special} are satisfied.
Moreover, the phase~${\mathfrak{r}}$ in~\cite[Section~3.4]{special} takes the values
\[ 3 {\mathfrak{r}} := \arg \left( -\frac{b}{t_0^3} \right) \approx 135^\circ \!\!\!\! \mod 180^\circ \:, \]
with an arbitrarily small error. Since~${\mathfrak{r}}$ must be chosen in the interval~$(-60^\circ, 0)$
(see~\cite[eqn~(3.35)]{special}), we conclude that
\[ {\mathfrak{r}} \approx -15^\circ \:. \]
Next, we consider the phase of the function~$h(z)$, which we write as
\[ h(z) = \pm 2\:\frac{e^{-2 i {\mathfrak{r}}}}{2^\frac{2}{3}}\: \frac{|t_0|^2}{t_0^2}\: |b|^{-\frac{2}{3}}
\:\sqrt{z^2-b}\: t_0\:. \]
It follows that
\[ \arg h(z) \approx -2 {\mathfrak{r}} \!\!\!\!\mod 180^\circ \approx 30^\circ \!\!\!\!\mod 180^\circ \:, \]
and thus~$\arg (h(z)^2) \approx 60^\circ$. This shows that the argument 
of the Airy function in~\cite[eqn~(3.36)]{special} does indeed avoid the branch cut.
\QED

\begin{Lemma} \label{lemma64}
For any~$c>0$, there are positive constants~$C_1$ and~$C_2$ such that
for sufficiently large~$|\Omega|$, the following statement holds. 
We consider the quadratic potential~\eqref{tVparabolic} with
parameters~${\mathfrak{p}}$, ${\mathfrak{q}}$ and~${\mathfrak{r}}$ in the range
\begin{align}
& \qquad\;\;\; |{\mathfrak{p}} {\mathfrak{q}}| \geq C_1\, |\Omega|^3 \label{c4} \\
|\re {\mathfrak{p}}| &\geq C_1\, |\Omega|\:, \qquad\quad\!\!
 |\im {\mathfrak{p}}| \leq c\, |\Omega| \label{c1} \\
|\re {\mathfrak{q}}| &\leq c\, |\Omega|^2\:, \qquad\quad
|\im {\mathfrak{q}}| \leq c\, |\Omega| \label{c2} \\
|\re {\mathfrak{r}}| &\leq c\:, \qquad \qquad \;\;\;\,
|\im {\mathfrak{r}}| \leq c\, |\Omega|^{-1} \:. \label{c3}
\end{align}
We choose~$\tilde{\phi}(u) = U_a^+(z(u))$ as the parabolic cylinder function
defined by the contour~$\Gamma_+ = \R + i$ (see~\cite[eqn~(3.2)]{special})
and let~$\tilde{y} = \tilde{\phi}'(u)/\tilde{\phi}(u)$ be the corresponding solution
of the Riccati equation. We denote the zero of~$\re \tilde{V}$ by~$u_1$ and set
\beq \label{upmdef}
u_\pm = u_1 \pm {C_2}\, |{\mathfrak{p}} {\mathfrak{q}}|^{-\frac{1}{6}}\:.
\eeq
Assume that~$z$ and~$b$ given by~\eqref{azdef} are in the range
\[ \im z^2, \im b > {C_2} \:. \]
Then for all~$u \in [0, \re {\mathfrak{r}}]$ we have the estimates
\begin{align}
|\re \tilde{y}| &\leq |\re \sqrt{\tilde{V}} | + C_1 \,|{\mathfrak{p}} {\mathfrak{q}}|^\frac{1}{6} \label{es1} \\
|\im \tilde{y}| &\leq \big| \im \sqrt{\tilde{V}} \big| + C_1\,|{\mathfrak{p}} {\mathfrak{q}}|^\frac{1}{6} \label{es3} \\
\frac{1}{2} \: \big| \im \sqrt{\tilde{V}} \big| &\leq -\im \tilde{y}
\qquad\text{if $u<u_-$} \:. \label{es2}
\end{align}
\end{Lemma}
\Proof Using the scaling of the parameters~${\mathfrak{p}}$, ${\mathfrak{q}}$ and~${\mathfrak{r}}$, we find
\begin{align*}
|u_1 - {\mathfrak{r}}| &\simeq \left| \frac{{\mathfrak{p}}}{{\mathfrak{q}}} \right|^\frac{1}{2} \\
\tilde{V}'(u_1) &= \frac{{\mathfrak{q}}}{2}\: (u_1 - {\mathfrak{r}}) \simeq |{\mathfrak{p}} {\mathfrak{q}}|^\frac{1}{2} \\
\tilde{V}'(u_1) \, (u_+-u_-) &\simeq C_2 \,|{\mathfrak{p}} {\mathfrak{q}}|^\frac{1}{3} \\
\tilde{V}''(u_1) \, (u_+-u_-)^2 &\simeq |{\mathfrak{q}}| \:C_2^2\: |{\mathfrak{p}} {\mathfrak{q}}|^{-\frac{1}{3}}
= C_2^2\: |{\mathfrak{p}} {\mathfrak{q}}|^\frac{1}{3} \left| \frac{{\mathfrak{q}}}{{\mathfrak{p}}^2} \right|^\frac{1}{3}
\end{align*}
In view of~\eqref{c1} and~\eqref{c2}, the quotient~${\mathfrak{q}}/{\mathfrak{p}}^2$ can be
made arbitrarily small by increasing~$C_1$. This makes it possible to arrange that
on the interval~$[u_-, u_+]$, the dominant term in~$\tilde{V}$ is the linear term. As a consequence,
\begin{align}
\tilde{V}(u_-) &\simeq \tilde{V}'(u_1)\, (u_- - u_1) \simeq {C_2}\, |{\mathfrak{p}} {\mathfrak{q}}|^\frac{1}{3}
\label{Vtes} \\
z^2 - b &= {\mathfrak{q}}^\frac{1}{2}\: (u-{\mathfrak{r}})^2 + 4 \frac{{\mathfrak{p}}}{\sqrt{{\mathfrak{q}}}} + 2
= \frac{4}{\sqrt{{\mathfrak{q}}}} \:\tilde{V} + 2 \\
z(u_-)^2 - b &\simeq {C_2}\, |{\mathfrak{p}}|^\frac{1}{3} \: |{\mathfrak{q}}|^{-\frac{1}{6}}
= {C_2}\, |b|^\frac{1}{3} \:.
\end{align}
Hence at~$u_-$, Lemma~\ref{lemmaaway} shows that the WKB approximation applies.
Possibly by increasing~${C_2}$, we can arrange that~$\tilde{y} = \pm \sqrt{\tilde{V}}$
with an arbitrarily small relative error. Clearly, this WKB estimate also holds for~$u<u_-$ and for~$u>u_+$.

In order to justify the sign in~\eqref{es2}, we choose the square roots such that
\[ \arg z \approx 45^\circ \:,\qquad \arg {\mathfrak{q}}^\frac{1}{4} \approx 135^\circ \:. \]
Then the WKB estimate of \cite[Theorem~3.3; see also eqn~(3.30)]{special} shows that
the function~$U_a^+$ is approximated by
\[ \tilde{U}_a(z) \sim \exp \left( \frac{\sqrt{z}}{4} \: \sqrt{z^2-b} \right) \:, \]
where the sign of the square root is chosen such that
\[ \arg \sqrt{z^2-b} \approx 45^\circ \quad \text{if~$u< u_-$} \qquad \text{and} \qquad
\arg \sqrt{z^2-b} \approx -45^\circ \quad \text{if~$u> u_+$} \:. \]
As a consequence,
\[ \tilde{y}(u) \simeq \frac{d}{du} \left( \frac{\sqrt{z}}{4} \: \sqrt{z^2-b} \right)
= \frac{1}{4}\: \frac{2z^2 - b}{\sqrt{z^2-b}}\: {\mathfrak{q}}^\frac{1}{4} \: .\]
A short calculation shows that
\[ \im \tilde{y}(u) < 0 \quad \text{if~$u< u_-$} \qquad \text{and} \qquad
\re \tilde{y}(u) > 0 \quad \text{if~$u> u_+$} \:. \]

It remains to estimate~$\tilde{y}$ on the interval~$[u_-, u_+]$. If this interval does not
intersect~$[0, \re {\mathfrak{r}}]$, there is nothing to do.
If this intersection is not empty and~$u_- \not \in [0, \re {\mathfrak{r}}]$, we replace~${\mathfrak{r}}$
by~${\mathfrak{r}}+1$. Thus we may assume that~$u_- \in [0, \re {\mathfrak{r}}]$.
In view of~\eqref{upmdef} and~\eqref{Vtes}, we know that
\[ \max_{[u_-, u_+]} \sqrt{|V|} \:(u_+-u_-) \simeq  {C_2}^\frac{3}{2} \:. \]
Hence we can apply Lemma~\ref{lemmagronwall} to~$\tilde{y}$ to obtain
\[ |\im \tilde{y}(u)| \geq \frac{1}{c_2}\: |\im \tilde{y}(u)| - c_2\: (u_+-u_-)\: \max_{[u_-,u_+]} |\im V| \]
with a constant~$c_2$ which depends only on~$C_2$.
From our assumption~\eqref{c1}--\eqref{c3} it follows that~$|\im V| \leq c' |\Omega|$, 
where~$c'$ depends only on~$c$. Moreover, at~$u_-$ we can use the WKB estimate
together with~\eqref{Vtes}. Also applying~\eqref{upmdef}, we obtain
\[ |\im \tilde{y}(u)| \geq \frac{1}{c_2^2}\: |{\mathfrak{p}} {\mathfrak{q}}|^{\frac{1}{6}}
\left( 1 - 2 \,c' \, c_2^2\: |\Omega|\: |{\mathfrak{p}} {\mathfrak{q}}|^{-\frac{1}{3}} \right) . \]
In view of~\eqref{c4}, by increasing~$C_1$ we can arrange that the
first summand dominates the second, meaning that
\[ |\im \tilde{y}(u)| \geq \frac{1}{2 c_2^2}\: |{\mathfrak{p}} {\mathfrak{q}}|^{\frac{1}{6}} \:. \]
Increasing~$C_1$ if necessary, we obtain the result.
\QED

We finally remark that it is a pure convention of the parabolic cylinder functions
defined in~\cite{special} that~$\tilde{y}$ lies in the lower half plane.
Solutions in the upper half plane are readily obtained with the following
double conjugation method: We consider the solution~$U^+_a(z)$ corresponding
to the complex conjugate potential~$\overline{\tilde{V}}$.
Then~$\phi(u):=\overline{U_a(z)^+}$ is a parabolic cylinder function corresponding to the
potential~$\tilde{V}$. The corresponding Riccati solution~$\tilde{y}(u) := \phi'(u)/\phi(u)$
satisfies~\eqref{es1}, \eqref{es3} and, in analogy to~\eqref{es2}, the inequality
\beq \frac{1}{2} \: \big| \im \sqrt{\tilde{V}} \big| \leq \im \tilde{y}
\qquad\text{if $u<u_-$} \:. \label{es2p}
\eeq

\subsection{Applications to the Angular Teukolsky Equation}
We now want to get estimates on an interval~$I=[u_{\min}, \umax]$ which
includes a zero of $\re V$, which we denote by~$u_1$.
We choose
\[ \tilde{V}(u) = V(u_1) + V'(u_1)\, (u-u_1) + \frac{1}{2}\: {\mathscr{V}}''\: (u-u_1)^2 \]
with
\[ {\mathscr{V}}'' := i \im V''(u_1) + \max_I \re V'' \:. \]
We define~$u_\pm$ as in Lemma~\ref{lemma64}.

\begin{Lemma} Assume that~$\im V \geq 0$ on the interval~$[u_{\min}, u_+]$
and that the assumptions of Lemma~\ref{lemma64} hold. Moreover, assume
that for sufficiently large~$|\Omega|$,
\beq \label{tilVlower}
|\tilde{V}(u)| \gtrsim |\Omega|^2\: (u-u_1)^2 \:.
\eeq
Then for sufficiently large~$|\Omega|$,
the invariant region estimate of Theorem~\ref{thmT} applies with~$g \equiv 0$ and
\[ \log T \big|^u_{u_{\min}} \leq C \, (u-u_{\min}) \qquad \text{for all~$u \in [u_{\min}, u_-]$}
\:, \]
where the constant~$C$ is independent of~$\Omega$.
\end{Lemma}
\Proof We take~$\tilde{V}$ as the approximate potential.
As the approximate solution~$\tilde{y}$ of the corresponding Riccati equation,
we take the the double conjugate solution introduced before~\eqref{es2p}.

The function~$f := \re(V-\tilde{V})$ has the properties
\[ f(u_1) = 0 = f'(u_1) \qquad \text{and} \qquad f''(u) \leq 0 \text{ on~$I$}\:. \]
Thus it is concave and lies below any tangent. In particular, it is everywhere negative,
\[ \re (V-\tilde{V}) \leq 0 \:. \]
Hence~\eqref{Urel} gives
\[ U \leq -\tilde{\beta}^2 \:. \]

We now estimate the error terms~$E_1, E_2, E_3$ in Theorem~\ref{thmT}:
\begin{align*}
\int_{u_{min}}^u |E_1| &\lesssim \int_{u_{min}}^u
\frac{1}{\tilde{\beta}(v)^2} \:\Big( |\alpha(v)|\: |\re(V-\tilde{V})| + |\re(V-\tilde{V})'| \Big) dv \\
&\lesssim \int_{u_{min}}^u
\frac{1}{\tilde{\beta}(v)^2} \:\Big( |\alpha(v)|\: |\Omega|^2\: |v-u_1|^3 + |\Omega|^2\: |v-u_1|^2 \Big) dv \\
&\lesssim |\Omega|^2 \int_{u_{min}}^u
\bigg( \frac{\alpha \tilde{\beta}}{\tilde{\beta}^3}\:  |v-u_1|^3 + \frac{|v-u_1|^2}{\tilde{\beta}^2} \bigg) dv \:. \\
\intertext{Applying Lemma~\ref{lemma64} gives}
&\lesssim |\Omega|^2 \int_{u_{min}}^u
\bigg( \frac{|\im \tilde{V}|}{|\tilde{V}|^\frac{3}{2}}\:  |v-u_1|^3 + \frac{|v-u_1|^2}{|\tilde{V}|} \bigg) dv \:. \\
\intertext{We now apply~\eqref{tilVlower} and use that~$|\im \tilde{V}| \lesssim |\Omega|$ to obtain}
&\lesssim |\Omega|^2 \int_{u_{min}}^u
\bigg( \frac{|\Omega|}{|\Omega|^3} + \frac{1}{|\Omega|^2} \bigg) dv \lesssim u-u_{\min}\:. \\
\intertext{Similarly,}
\int_{u_{min}}^u |E_2| &\lesssim \int_{u_{min}}^u \frac{1}{|\tilde{\beta}|}\:|\im(V-\tilde{V})|\:dv
\lesssim \int_{u_{min}}^u \frac{1}{|\tilde{\beta}|}\: |\Omega|\, |v-u_1|^3\:dv \\
&\lesssim \int_{u_{min}}^u |v-u_1|^2\:dv \lesssim u-u_{\min} \\
\int_{u_{min}}^u |E_3| &\lesssim \int_{u_{min}}^u \frac{|\im V|}{|\tilde{\beta}|^3}\: |\re(V-\tilde{V})|\: dv \\
&\lesssim |\Omega|^3 \int_{u_{min}}^u \frac{|v-u_1|^3}{|\tilde{V}|^\frac{3}{2}}\: dv 
\lesssim u-u_{\min}\:.
\end{align*}
This completes the proof.
\QED

\section{Estimates for a Singular Potential} \label{sec7}
At~$u=0$, the potential~\eqref{Vdef} has a pole of the form
\[ V(u) = \left( (k-s)^2 - \frac{1}{4} \right) \frac{1}{u^2} + \O(u^0) \:. \]
In preparation for estimating the solutions near this pole (see Section~\ref{secpoles} below),
in this section we analyze solutions of the Riccati equation for a potential includes the pole
and involves a general constant. More precisely, setting~$L=|k-s|$, we consider a potential of the form
\beq \label{Vtpole}
V(u) =  \left( L^2 - \frac{1}{4} \right) \frac{1}{u^2} + \zeta^2
\eeq
for a complex parameter~$\zeta$ and a non-negative integer~$L$.
In the case~$L=0$, the real part of~$V$ tends to~$-\infty$ as~$u \searrow 0$, whereas in the case~$L>0$,
it tends to~$+\infty$. We treat these two cases separately.

\subsection{The case~$L=0$}
In this case, the potential~\eqref{Vtpole} becomes
\beq \label{VL0}
V(u) = - \frac{1}{4 u^2} + \zeta^2 \:.
\eeq
We assume that~$\zeta$ lies in the upper right half plane excluding the real axis,
\[ \arg \zeta \in (0, 90^\circ] \:. \]
The corresponding Sturm-Liouville equation~\eqref{SLtilde} has explicit solutions
in terms of the Bessel function~$K_0$ and~$I_0$ (see~\cite[\S10.2.5]{DLMF}).
We choose
\beq \label{K0}
\phi(u) = -\sqrt{u} \:\Big( K_0(\zeta u) + \big( \arg \zeta - \log(2) + \gamma + i \big) \,I_0(\zeta u) \Big)
\eeq
(where~$\gamma \approx 0.577$ is Euler's constant).
Near the origin, we have the asymptotics~\cite[eq.~(10.31.2)]{DLMF}
\[ \phi(u) = \sqrt{u}\: \log |\zeta u| + i \sqrt{u} + \O(u) \:. \]
For large~$u$, on the other hand, we have the asymptotics (see~\cite[\S10.40(i)]{DLMF})
\[ \phi(u) = -\frac{e^{\zeta u}}{\sqrt{2 \pi \zeta}} \:\big( 
\arg \zeta - \log(2) + \gamma + i \big)\:\big(1+\O\big( (\zeta u)^{-1} \big) \big) . \]
We again denote the corresponding solution of the Riccati equation by~$y=\phi'/\phi$.

\begin{Prp} \label{prp71}
On the interval
\beq \label{usmall}
0 \leq u \leq \frac{1}{2 \,|\zeta|} \:,
\eeq
the $T$-method of Theorems~\ref{thm2} and~\ref{thmT} applies with~$g \equiv 0$ and
\beq \label{Vtphiform}
\tilde{V} = -\frac{1}{4u} \:, \qquad
\tilde{\phi}(u) = \sqrt{u}\: \log |\zeta u| + i \sqrt{u} \:.
\eeq
Moreover, the function~$T$ is bounded uniformly in~$\zeta$.

If in addition~$\im \zeta^2 > 0$, there is a constant~$C$ which
depends only on~$\arg \zeta$ such that
\beq \label{Imylower}
|\re y| \leq C \,|\zeta| \qquad \text{and} \qquad \im y \geq \frac{|\zeta|}{C} \qquad \text{for
all~$u > \frac{1}{2 \,|\zeta|}$}\:.
\eeq
\end{Prp}
\Proof Introducing the rescaled variable~$u' = |\zeta|\, u$, one sees that it suffices
to consider the case~$|\zeta|=1$. Choosing~$\tilde{V}$ and~$\tilde{\phi}$
on the interval~$(0, \frac{1}{2}]$ as in~\eqref{Vtphiform}, we obtain
\begin{align*}
\tilde{y}(u) &= \frac{i-2-\log u}{2 u \,(i-\log u)} \\
\alpha(u) &= \frac{(1+ \log u)^2}{2 u \,(1 + \log^2 u)} \leq \frac{1}{u} \:,\qquad
\tilde{\beta}(u) = \frac{1}{u \,(1+\log^2 u)} \\
U &= \re(V-\tilde{V}) - \tilde{\beta}^2 = \re \zeta^2 - \frac{1}{u^2 \,(1+\log^2 u)^2} \:. \\
\intertext{Since on the interval~$(0, \frac{1}{2}]$, the inequality~$2 u^2 \,(1+\log^2 u)^2<1/5$ holds,
we conclude that~$U<0$ and}
|U| &\geq \frac{1}{5 u^2 \,(1+\log^2 u)^2} \\
|E_1| &\leq 10\, u \,(1+\log^2 u)^2 \: \big| \re (\zeta^2) \big| \\
|E_2| &\leq 5 u \,(1+\log^2 u)\: \big| \im (\zeta^2) \big| \\
|E_3| &\leq 5 u^3 \,(1+\log^2 u)^3\: \big| \re(\zeta^2)\: \im(\zeta^2) \big| \:.
\end{align*}
This shows that Theorem~\ref{thmT} applies and that the function~$T$ is uniformly bounded.

If~$\im (\zeta^2)>0$, Lemma~\ref{lemma36} shows that the solution stays in the upper half plane
(this can also be seen directly from the differential equation~\eqref{riccati}).
Moreover, we know that at~$u=1/2$, the function~$y$ is bounded.
Furthermore, in the limit~$u \rightarrow \infty$, the solution~$y$ tends to the
stable fixed point~$\zeta$ (where we choose the sign of~$\zeta$ such that~$\re \zeta>0$;
for details see~\cite[Section~2]{invariant}).
Hence there is~$C>0$ such that~$|\re y| \leq C$
and~$\im y \geq 1/C$ on~$[1/2, \infty)$. This concludes the proof.
\QED

\subsection{The case~$L>0$}
We now consider the potential~\eqref{Vtpole} in the case~$L>0$. We assume that~$\zeta$
does not lie on the positive real axis,
\[ \arg \zeta \not \in 0 \!\!\! \mod 2 \pi \:. \]
The corresponding Sturm-Liouville equation~\eqref{SLtilde} has an explicit solution
in terms of the Bessel function~$K_L$ (see~\cite[\S10.2.5]{DLMF}), 
\beq \label{phibessel}
\phi(u) = \sqrt{u}\, K_L(-\zeta u) \:.
\eeq
Using the recurrence relations in~\cite[eqs.~(10.29.2) and~(10.29.1)]{DLMF}, it follows that
\beq \label{ybessel}
y(u) = \frac{1-2L}{u} + \zeta\: \frac{K_{L-1}(-\zeta u)}{K_{L}(-\zeta u)} \:.
\eeq

Near the origin, we have the asymptotics (see~\cite[eqs.~(10.31.1) and~(10.25.2)]{DLMF})
\[ \phi(u) = \frac{(n-1)!}{\sqrt{2 \zeta}} \left( -\frac{\zeta u}{2} \right)^{\frac{1}{2}-L} (1 + \O(u))
\;-\; \frac{\sqrt{2}}{L!\, \sqrt{\zeta}} \left( \frac{\zeta u}{2} \right)^{\frac{1}{2}+L} (1 + \O(u)) \:, \]
whereas for large~$u$, we have the asymptotics (see~\cite[eq.~(10.25.3)]{DLMF})
\[ \phi(u) \sim \sqrt{-\frac{\pi}{2 \zeta}}\: e^{\zeta u}  \:,\quad
y(u) \sim \zeta \qquad \text{as~$u \rightarrow +\infty$} \]
(where the square root is taken such that~$\re \sqrt{-\zeta} > 0$).

\begin{figure}
\begin{center}
\includegraphics[width=9cm]{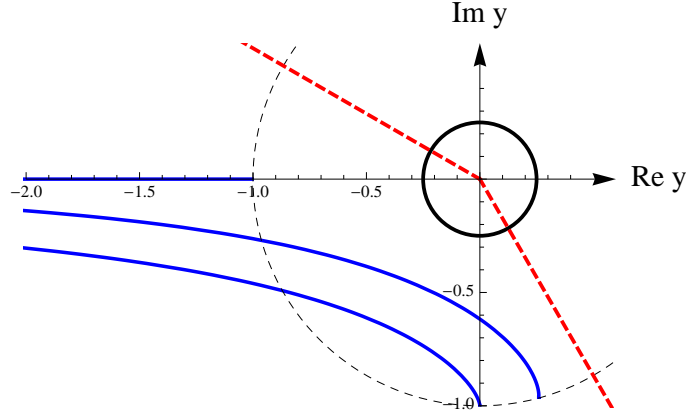}
\caption{The Bessel solution~$y$ for~$\zeta=-1$, $\zeta=-i$
and~$\zeta = \exp(-75\, i \pi /180)$.}
\label{figzeta}
\end{center}
\end{figure}

\begin{Prp} \label{prp72} For sufficiently small~$\varepsilon>0$, the following statement holds:
If the argument of~$\zeta$ lies in the range
\beq \label{argcond}
\arg \zeta \in \big( 180^\circ-\varepsilon, 180^\circ+\varepsilon \big) \cup
\big(-90^\circ-\varepsilon, 90^\circ+\varepsilon \big) \cup \big( -75^\circ-\varepsilon, 75^\circ+\varepsilon \big) \:,
\eeq
then the Bessel solution~\eqref{ybessel} satisfies for all~$u \in \R^+$ the inequalities
\beq \label{conds}
|y(u)| \geq \frac{|\zeta|}{4} \qquad \text{and} \qquad
150^\circ < \arg y(u) < 300^\circ \:.
\eeq
\end{Prp}
\Proof Rescaling the variables by~$|\zeta|\, u \rightarrow u$,
we may assume that~$|\zeta|=1$. The solutions for~$\zeta=-1$, $\zeta=-i$
and~$\zeta = \exp(-75\, i \pi /180)$ satisfy~\eqref{conds}, as one sees
from Figure~\ref{figzeta}.
The result now follows by continuity.
\QED

\section{Estimates for the Angular Teukolsky Equation near the Poles} \label{secpoles}
Near~$u=0$, the potential~\eqref{Vdef} has the expansion
\[ V(u) = \left( L^2 - \frac{1}{4} \right) \frac{1}{u^2} + c_0 + c_2 u^2 + \O(\Omega^2 u^4) \:, \]
where we again set~$L=|k-s|$, and where the coefficients~$c_0$ and~$c_2$ scale in~$\Omega$ like
\beq \label{Oscale}
c_0 = -2 s \Omega - \mu + \O(\Omega^0)\:, \qquad c_2 = \Omega^2 + \O(\Omega) \:.
\eeq
We again treat the cases~$L=0$ and~$L>0$ separately.

\subsection{The Case~$L=0$}
Our goal is to estimate the solutions on an interval~$I:=(0, \umax]$.
We choose the approximate potential according to~\eqref{VL0},
\[ \tilde{V}(u) = - \frac{1}{4 u^2} + \zeta^2 \:, \]
and take~\eqref{K0} as the solution~$\tilde{\phi}$ of the corresponding
Sturm-Liouville equation~\eqref{SLtilde}.
The constant~$\zeta$ in~\eqref{Vtpole} is chosen as
\[ \im \zeta^2 = \im c_0 \qquad \text{and} \qquad
\re \zeta^2 = \max_{I} \re \left( V(u) -  \Big( L^2 - \frac{1}{4} \Big) \frac{1}{u^2} \right) . \]
Then by construction we have~$\re(V-\tilde{V}) \leq 0$.
In view of~\eqref{Urel}, we conclude that~$U$ is negative, making it possible to
apply the $T$-method. Moreover,
\begin{align}
|U| &\geq \tilde{\beta}^2\:, & |U| &\geq 2 \tilde{\beta}\: \sqrt{\re(V-\tilde{V})} \label{Ues} \\
|\re (V-\tilde{V})| &\lesssim |\Omega|^2 \:\umax^2 \:,& |\re (V-\tilde{V})'| &\lesssim |\Omega|^2 \: \umax
\label{VVtes} \\
|\im (V-\tilde{V})| &\lesssim |\Omega|\: u^2 \label{ImVes}
\end{align}

\begin{Prp} Assume that~$\im \zeta^2 > 0$. We consider the solution~$\phi$
on the interval
\beq \label{umaxdef}
(0, \umax] \qquad \text{with} \qquad \umax \lesssim |\Omega|^{-\frac{1}{2}} \:,
\eeq
having the following asymptotics near~$u=0$,
\[ \phi(u) = \sqrt{u}\: \log |\zeta u| + i \sqrt{u} + \O(u) \:. \]
Then the $T$-estimates of Theorems~\ref{thm2} and~\ref{thmT} apply with~$g \equiv 0$ and
\[ \log T(u) \leq C \, |\Omega| \, \umax^2 \left(1 + \log^4 |\zeta u| \right) 
\qquad \text{for all~$u \in (0, \umax]$}\:. \]
Here
\[ \left\{ \begin{array}{cl} \text{$C$ is a numerical constant} & \text{if~$ \umax \leq (2|\zeta|)^{-1}$} \\[0.2em]
\text{$C$ depends on~$\arg \zeta$} & \text{if~$\umax > (2|\zeta|)^{-1}$\:.} 
\end{array} \right. \]
\end{Prp}
\Proof The function~$\tilde{\phi}$ coincides with the function~$\phi$ in
Proposition~\ref{prp71}.
Hence on the interval~\eqref{usmall}, we obtain, for a suitable numerical constant~$c$,
\begin{align*}
|\alpha| &\leq \frac{c}{u}\:, \qquad |\tilde{\beta}| \geq \frac{1}{c u \,(1+\log^2 \!|\zeta u|)} \\
|E_1| &\leq \frac{2 c^3}{u}\: u^2 \,(1+\log^2 \!|\zeta u|)^2\: |\re (V-\tilde{V})|
+ \frac{ c^2 u^2}{2}\: (1+\log^2 \!|\zeta u|)^2\: |\re (V-\tilde{V})'| \\
&\overset{\eqref{VVtes}}{\lesssim} \Big( 2 c^3 + \frac{ c^2}{2} \Big)\:
|\Omega|^2\: \umax^3\:(1+\log^2 \!|\zeta u|)^2 \\
|E_2| &\leq \frac{1}{\tilde{\beta}}\, |\im (V-\tilde{V})|
\overset{\eqref{ImVes}}{\lesssim} c \,|\Omega|\: u^3\:(1+\log^2 \!|\zeta u|) \\
|E_3| &\leq \frac{|\im V|}{\tilde{\beta}^3}\: |\re(V-\tilde{V})| 
\overset{\eqref{VVtes}}{\lesssim} c^3\:|\Omega|^3\: u^3 \umax^2\,(1+\log^2 \!|\zeta u|)^2 \:.
\end{align*}

In the remaining region~$u>1/(2 |\zeta|)$, we know in view of~\eqref{umaxdef} that
\beq \label{zetalower}
|\zeta| \gtrsim |\Omega|^\frac{1}{2}\:.
\eeq
We have the global estimates~\eqref{Imylower} for~$\alpha$ and~$\tilde{\beta}$.
Using the second inequality in~\eqref{Ues}, we can estimate the error terms by
\begin{align*}
|E_1| &\leq \frac{\alpha}{\tilde{\beta}}\: \sqrt{\re (V-\tilde{V})}
+ \frac{|\re (V-\tilde{V})'|}{\tilde{\beta}^2}  \\
&\lesssim C^2 \,|\Omega|\, \umax + \frac{C^2}{|\zeta|^2}\: |\Omega|^2 \, \umax 
\overset{\eqref{zetalower}}{\lesssim}
C^2 \,|\Omega|\, \umax \\
|E_2| &\leq \frac{|\im (V-\tilde{V})|}{|\tilde{\beta}|} \lesssim \frac{C |\Omega|}{|\zeta|} 
\overset{\eqref{zetalower}}{\lesssim} C\, |\Omega|^{\frac{1}{2}} \\
|E_3| &\leq \frac{|\im V|}{\tilde{\beta}^3}\: |\re(V-\tilde{V})| \lesssim
\frac{C^3}{|\zeta|^3}\: |\Omega|^3\, \umax^2
\overset{\eqref{zetalower}}{\lesssim} 
C^3\: |\Omega|^{\frac{3}{2}}\, \umax^2 \:.
\end{align*}

Integrating the error terms from~$u$ to~$\umax$, using~\eqref{umaxdef}
and renaming the constants gives the result.
\QED

\subsection{The Case~$L>0$}
We choose~$u_0$ such that
\[ \re V'(u_0) = 0 \:. \]
According to~\eqref{Oscale}, we know that
\beq \label{u0bound}
u_0 \simeq |\Omega|^{-\frac{1}{2}} \:.
\eeq
Our goal is to estimate the solutions on the interval~$I:=(0, u_0]$.
We take~\eqref{Vtpole} as our approximate potential
\beq \label{tilV}
\tilde{V} = \left( L^2 - \frac{1}{4} \right) \frac{1}{u^2} + \zeta^2 \:,
\eeq
where
\beq \label{zetarel}
\zeta^2 = c_0 - (1+2 i) \, {\mathcal{C}}^2 \,|\Omega| \:,
\eeq
and~${\mathcal{C}}$ is a constant to be determined later.
Then
\beq \label{VVtansatz}
V - \tilde{V} = (1+2i) \,{\mathcal{C}}^2 \,|\Omega| + c_2 u^2 + \O(\Omega^2 u^4) \:.
\eeq
We choose a constant
\beq \label{c3choice}
c_3 \geq \frac{1}{|\Omega|} \sup_{(0, u_0]} \left| V - \tilde{V} - (1+2i) \,{\mathcal{C}}^2 \,|\Omega| \right| \:.
\eeq
In view of~\eqref{Oscale} and~\eqref{u0bound}, the constant~$c_3$ can indeed be chosen
independent of~$\Omega$.

\begin{Lemma} \label{lemmazeta}
For every~${\mathcal{C}}_\text{min}>0$ there is~${\mathcal{C}}_\text{max}$
such that for every~$c_0$ with
\beq \label{c0bound}
|\im c_0| \lesssim |\Omega| \:,
\eeq
there is a parameter~${\mathcal{C}}$ in the range
\[ {\mathcal{C}}_\text{min} \leq {\mathcal{C}} \leq {\mathcal{C}}_\text{max} \]
and a complex number~$\zeta$ which satisfies~\eqref{zetarel}, such that
the conditions in Proposition~\eqref{prp72} hold.
\end{Lemma}
\Proof In the limit~$|c_0| \rightarrow \infty$, the real part of~$c_0$ dominates
its imaginary part in view of~\eqref{c0bound}. Thus taking~${\mathcal{C}}={\mathcal{C}}_\text{min}$,
the argument of~$\zeta^2$ tends to zero or~$\pi$ as~$|c_0| \rightarrow \infty$.
Taking the square root, we can thus satisfy~\eqref{argcond}.
More precisely, there is a constant~$c_4$ (independent of~$\Omega$)
such that~\eqref{argcond} holds if~$|c_0| > c_4\, |\Omega|$.

In the remaining case~$|c_0| \leq c_4\, |\Omega|$, we choose~${\mathcal{C}} = {\mathcal{C}}_\text{max}$.
By choosing~${\mathcal{C}}_\text{max} \gg c_4$, we can arrange that the argument of~$\zeta^2$
lies arbitrarily close to~$\arg(-(1+2i)$. Taking the square root, we see that~\eqref{argcond} again holds.
\QED

We are now in a position to apply the $\kappa$-method of Proposition~\ref{prpnewkappa}.
We choose the approximate potential~\eqref{tilV} with~$\zeta$
in agreement with Lemma~\ref{lemmazeta}.
Moreover, we choose the solution~$\tilde{\phi}$ of the corresponding
Sturm-Liouville equation~\eqref{SLtilde} to be the Bessel solution~\eqref{phibessel}.
Again, we denote the corresponding Riccati solution by~$\tilde{y}=\tilde{\phi}'/\tilde{\phi}$.
It has the properties
\beq \label{ytprop}
|\tilde{y}| \geq \frac{{\mathcal{C}}}{4}\: |\Omega|^\frac{1}{2} 
\eeq
and
\beq \label{abprop}
150^\circ < \arg \tilde{y} < 300^\circ \:.
\eeq

\begin{Prp} Choosing the function~$g$ in~\eqref{kset} as
\beq \label{gdef2}
g(u) = {\mathcal{C}}^2 \, |\Omega|^\frac{1}{2}\: \sigma(u) \:,
\eeq
and choosing~${\mathcal{C}}$ sufficiently large (independent of~$\Omega$),
the disks defined by~\eqref{ninvdisk} are invariant under the backward
Riccati flow. Moreover, the determinator~${\mathfrak{D}}$ is positive.
\end{Prp}
\Proof We want to apply Proposition~\ref{prpnewkappa}
starting at~$u_0$ going backwards to the singularity at the origin.
Up to now, we always applied the invariant region estimates for increasing~$u$.
In order to avoid confusion, we now replace~$u$ by~$-u$, so that we need to estimate the solution on
the interval~$[-u_0, 0)$ (note that the potential in the Sturm-Liouville equation
does not change sign under the change of variables~$u \rightarrow -u$).
Since the relation~$\tilde{y}=\tilde{\phi}'/\tilde{\phi}$ involves one derivative,
the function~$\tilde{y}$ changes sign, so that~\eqref{abprop} becomes
\beq \label{abprop2}
-30^\circ \leq \arg \tilde{y} \leq 120^\circ \:.
\eeq

Using~\eqref{VVtansatz} and~\eqref{c3choice}, we obtain
\begin{align}
|\re(V-\tilde{V})'| &\lesssim |\Omega|^2 |u|  \lesssim |\Omega|^\frac{3}{2} \label{VVtprime} \\
2 \alpha \re (V-\tilde{V}) &\geq 2 \alpha \:{\mathcal{C}}^2 \,|\Omega| - 2 \alpha\, c_3\, |\Omega| \nonumber \\
\tilde{\beta} \im (V-\tilde{V}) &\geq 2 \tilde{\beta} \:{\mathcal{C}}^2 \,|\Omega| - \tilde{\beta} \,c_3\, |\Omega| \nonumber \:.
\end{align}
Adding the last two inequalities and using the trigonometric bound
\[ \alpha + \tilde{\beta} = \left\langle \begin{pmatrix} 1 \\ 1 \end{pmatrix},
\begin{pmatrix} \alpha \\ \tilde{\beta} \end{pmatrix} \right\rangle 
\overset{\eqref{abprop2}}{\geq} \sqrt{2}\, |\tilde{y}|\, \cos 75^\circ \geq \frac{1}{3}\, |\tilde{y}| \:, \]
we get
\begin{align*}
2 \alpha & \re (V-\tilde{V})+\tilde{\beta} \im (V-\tilde{V}) \geq
\frac{2}{3}\, |\tilde{y}|\, {\mathcal{C}}^2 \,|\Omega| - 3\, |\tilde{y}| \,c_3\, |\Omega| \\
& = \frac{2}{3}\, |\tilde{y}|\,( {\mathcal{C}}^2 - 5 c_3) \,|\Omega| 
\overset{\eqref{ytprop}}{\geq} \frac{\mathcal{C}}{6}\: ({\mathcal{C}}^2 - 5 c_3) \,|\Omega|^\frac{3}{2} \:.
\end{align*}
By choosing~${\mathcal{C}}$ sufficiently large (independent of~$|\Omega|$), we can arrange that
\beq \label{arrange}
2 \alpha \re (V-\tilde{V})+\tilde{\beta} \im (V-\tilde{V}) \geq \frac{{\mathcal{C}}^3}{8}\: |\Omega|^\frac{3}{2} \:.
\eeq

Note that the function~$g$ in~\eqref{gdef2} is monotone increasing because of~\eqref{sigmadef}
and~\eqref{abprop2}. Moreover,
\[ \kappa = {\mathcal{C}}^2 \, |\Omega|^\frac{1}{2} + 
\frac{1}{\sigma} \int_{u_0}^u \sigma\: \im(V-\tilde{V}) \:. \]
Using~\eqref{VVtansatz} and~\eqref{c3choice} together with the fact that~$\sigma$ is
monotone increasing, we conclude that
\[  {\mathcal{C}}^2 \,|\Omega|^\frac{1}{2} - c_3 \big| \Omega u_0 \big| \leq \kappa \leq {\mathcal{C}}^2 \big(|\Omega|^\frac{1}{2} + |\Omega u_0| \big) + c_3 \big|\Omega u_0 \big| \:. \]
Using~\eqref{u0bound}, possibly by increasing~${\mathcal{C}}$ we can arrange that
\[ \frac{{\mathcal{C}}^2}{2}\:|\Omega|^\frac{1}{2} \leq \kappa \leq 2 \,{\mathcal{C}}^2\, |\Omega|^\frac{1}{2} \:. \]
We now estimate~$\kappa - R$ using Lemma~\ref{Rpmb}. Keeping in mind that~$\tilde{\beta} \geq 0$,
we obtain
\begin{align*}
\frac{|\re(V-\tilde{V})|}{|\tilde{\beta} + \kappa|} &\leq
\frac{2 ({\mathcal{C}}^2 +c_3)}{{\mathcal{C}}^2} \, |\Omega|^\frac{1}{2} \\
\frac{\kappa^2}{|\tilde{\beta} + \kappa|} &\leq \kappa \leq 2 \,{\mathcal{C}}^2\, |\Omega|^\frac{1}{2} \:.
\end{align*}
Thus, possibly after increasing~${\mathcal{C}}$, we obtain
\[ |\kappa-R| \leq 2 \,{\mathcal{C}}^2\, |\Omega|^\frac{1}{2} \:. \]
As a consequence,
\beq \label{kmR}
|(\kappa-R) \im V| \lesssim 2 \,{\mathcal{C}}^2\, |\Omega|^\frac{3}{2} \:.
\eeq

Comparing~\eqref{arrange} with~\eqref{VVtprime} and~\eqref{kmR}, one sees
that, possibly by further increasing~${\mathcal{C}}$, we can arrange 
that the determinator as given by~\eqref{Ddef3} is positive.
This concludes the proof.
\QED

\Thanks {{\em{Acknowledgments:}}
We are grateful to the Vielberth Foundation, Regensburg, for generous support.

\providecommand{\bysame}{\leavevmode\hbox to3em{\hrulefill}\thinspace}
\providecommand{\MR}{\relax\ifhmode\unskip\space\fi MR }
\providecommand{\MRhref}[2]{%
  \href{http://www.ams.org/mathscinet-getitem?mr=#1}{#2}
}
\providecommand{\href}[2]{#2}

\end{document}